\documentclass[11pt]{article}
\usepackage{amssymb}
\usepackage{graphicx}
\usepackage{setspace}
  \usepackage{paralist}
  \usepackage{longtable}
   \usepackage{multirow}
    \usepackage{rotating}
\usepackage{fancyhdr}

\usepackage{amsmath, amsthm, amssymb}
\usepackage{authblk}
\usepackage{graphicx}
\usepackage{float}
\usepackage{hyperref}
\usepackage[margin=1in]{geometry}
\usepackage{comment}
\usepackage{color}
\allowdisplaybreaks[4]
\numberwithin{equation}{section}
\date{}

\newtheorem{theorem}{Theorem}[section]

\newtheorem{lemma}[theorem]{Lemma}

\newtheorem{remark}[theorem]{Remark}

\newcommand{\be}{\begin{equation}}
\newcommand\ee{\end{equation}}
\newcommand\bes{\begin{eqnarray}}
\newcommand\ees{\end{eqnarray}}
\newcommand\bess{\begin{eqnarray*}}
\newcommand\eess{\end{eqnarray*}}
\newcommand\D{\displaystyle}

\title{Pointwise space-time estimates of two-phase fluid model in dimension three}

\author{Zhigang Wu\thanks{Corresponding to:\ zgwu@dhu.edu.cn},\ \ Wenyue Zhou
\vspace{1mm}\\
\textit{Department of Mathematics, Donghua University,}\\\textit{\small Shanghai 201620, P.R. China}
}

\begin{document}

\pagestyle{myheadings} \markboth{WU & Zhou}{Two-Phase Fluid Model}\maketitle
\renewcommand{\thefootnote}{\fnsymbol{footnote}}


\renewcommand{\thefootnote}{\arabic{footnote}}


{{\bf  Abstract:}} In this paper, we investigate the pointwise space-time behavior of  two-phase fluid model derived by Choi \cite{Choi} [SIAM J. Math. Anal., 48(2016), pp. 3090-3122], which is the compressible damped Euler equations coupled with compressible Naiver-Stokes equations. Based on Green's function method together with frequency analysis and nonlinear coupling of different wave patterns, it shows that both of two densities and momentums obey the generalized Huygens' principle as the compressible Navier-Stokes equations \cite{LW}, however, it is different from the compressible damped Euler equations \cite{Wang2}. The main contributions include seeking suitable combinations to avoid the singularity from the Hodge decomposition in the low frequency part of the Green's function, overcoming the difficulty of the non-conservation arising from the damped mechanism of the system, and developing the detailed description of the singularities in the high frequency part of the Green's function. Finally, as a byproduct, we extend $L^2$-estimate in \cite{Wugc} [SIAM J. Math. Anal., 52(2020), pp. 5748-5774] to $L^p$-estimate with $p>1$.



\section{Introduction}

In this paper, we investigate the two-phase fluid model consisting of the compressible isothermal Euler equations coupled with compressible isentropic Naiver-Stokes equations, which reads as
\begin{equation}\label{1.1}
\left\{\begin{array}{l}
\rho_{t}+\operatorname{div}(\rho u)=0, \\
(\rho u)_{t}+\operatorname{div}(\rho u \otimes u)+\nabla \rho=-\rho(u-v), \\
n_{t}+\operatorname{div}(n v)=0,\\
(n v)_{t}+\operatorname{div}(n v \otimes v)+\nabla P(n)-\mu \Delta v-(\mu+\lambda) \nabla \operatorname{div} v=\rho(u-v),\ \ (x,t)\in \mathbb{R}^3\times\mathbb{R}^+.
\end{array}\right.
\end{equation}
Here the unknowns $\rho(x,t)$, $n(x,t)$ are the densities of fluid, and $u(x,t)$, $v(x,t)$ are their corresponding velocities of $\rho (x,t)$ and $n(x,t)$. The pressure function $P(n)=An^\gamma$ ($A>0, \gamma\geq1$). The parameters $\mu$ and $\lambda$ are the shear viscosity coefficient and the bulk viscosity coefficient respectively, and satisfy the physical conditions: $\mu>0$ and $\frac{2}{3}\mu+\lambda\geq0$.

The model (\ref{1.1}) was first derived  by taking the hydrodynamic limit from the Vlasov-Fokker-Planck/isentropic Navier-Stokes equations (VFPNS) with local alignment forces in \cite{Choi}. VFPNS system in \cite{Choi} is a kinetic-fluid model and used to motion the interactions between particles and fluid, which has been attracted more attention for applications in dynamic spray, diesel engines, etc.\cite{Berres,Sartory,Tory,Williams}. For VFPNS system without the local alignment forces, Mellet and Vasseur \cite{Mellet} studied the global weak solution in bounded domains. In terms of the other related model, we refer to Baranger and Desvillettes \cite{Baranger} for the local existence of the classical solution to the Vlasov/compressible Euler equations, Duan and Liu \cite{Duan} for the global existence and decay rate of the classical solution with the small initial data to the Vlasov-Fokker-Planck/compressible Euler equations,  Carrillo $et\ al.$ \cite{Carrillo} for  a coupled kinetic-fluid model consisting of the isothermal Euler equations and incompressible Navier-Stokes equations, Choi \cite{Choi2}, Karper $et\ al.$ \cite{Karper1} for the global classical and weak solutions of the Vlasov-Fokker-Planck equation, respectively.

There are a lot of works on the global existence and large time behavior of the classical solution with small initial perturbation for the Cauchy problem of the compressible fluid models. We just review two of the most closely related to the model (\ref{1.1}). For the compressible Navier-Stokes equations in 3D, the global existence and $L^2$-decay rate of the solution were given in Matsumura and Nishida \cite{Matsumura1,Matsumura2}, and the $L^p$-decay rate with $p\geq2$ was established  in Ponce \cite{Ponce}. When there is a potential force term, the similar results was developed by Duan $et\ al.$ in \cite{Duan2,Duan3}. Later on, Li and Zhang \cite{Lihl} improved the decay rate when the initial data belongs to some suitable negative Sobolev space, Guo and Wang \cite{Guo} directly obtained $L^2$-decay rate by a pure energy method without the spectrum analysis. For the compressible Euler equation with damping, Wang and Yang \cite{Wang2} investigated the global existence and the pointwise estimates in $H^4$-framework by using the Green's function method together with the energy method, Sideris $et\ al.$ in \cite{Sideris}, Chen and Tan \cite{Chen}, Tan and Wu \cite{Tan} used the different methods to obtain the global existence and $L^2$-decay rate in $H^3$-framework.

For the Cauchy problem of the system (\ref{1.1}), there are few results. The first result was given in Chio \cite{Choi} for the small initial perturbation in $H^l$-space, where they obtained the global existence and the $L^2$-decay rate for the periodic domain, and the global existence in the whole space. Later on, Wu $et\ al.$ \cite{Wugc} resolved the whole space problem and gave the optimal $L^2$-decay rate of the solution and its higher-order spatial derivatives, where they used Hodge decomposition and the spectral analysis. Very recently, Tang and Zhang \cite{Tang} reconsidered the Cauchy problem by the spectral analysis but without using the Hodge decomposition.

In order to explicitly describe the wave propagation of the compressible fluid model, one needs to study the space-time pointwise estimates, since the usual $L^2$-estimates only exhibit the dissipative properties of solutions through the energy method combining with the spectrum analysis. The pioneering works in this direction are Zeng \cite{Z1} and Liu and Zeng \cite{LZ2} for 1D compressible fluid models. The isentropic compressible Navier-Stokes system in 3D was investigated by Hoff and Zumbrun \cite{HZ0,HZ} and Liu and Wang \cite{LW} for the linear and nonlinear problems respectively.  A hyperbolic-parabolic system obeys the generalized Huygens' principle in \cite{LW} implies that its pointwise space-time description of the solution contains both a diffusion wave (D-wave): $(1+t)^{-\frac{3}{2}}\big(1+\frac{|x|^2}{1+t}\big)^{-\frac{3}{2}}$ and a generalized Huygens' wave (H-wave): $(1+t)^{-2}\big(1+\frac{(|x|-ct)^2}{1+t}\big)^{-\frac{3}{2}}$ in $\mathbb{R}^3$. Obviously, the $L^2$-decay rate of these two waves is the same as the heat kernel, and D-wave decays faster than H-wave when $p\leq2$, H-wave decays faster than D-wave when $p\geq2$  for the $L^p(\mathbb{R}^3)$-estimates. On the other hand, due to the damped mechanism, the pointwise estimates for the damped Euler equation in \cite{Deng2,Wang2} does not obey the generalized Huygens' principle, which is different from the compressible Navier-Stokes equations.  Later on, there also are series of results in this direction based on the Green's function method for other compressible fluid models, for instance, the non-isentropic Navier-Stokes equations in \cite{Du,LD}, the unipolar Navier-Stokes-Poisson equations in \cite{Wang} and references therein. These results imply that the different models exhibit different wave propagations or wave patterns.

The goal of this paper is to derive the generalized Huygens' principle for the system \eqref{1.1}. As we know, when deducing the generalized Huygens' principle for compressible fluid models, such as the Navier-Stokes equations \cite{Deng,LS,LW},  the bipolar Navier-Stokes-Poisson equations \cite{Wu4,Wu7}, the micropolar fluid model \cite{Wu6}, the conservation is critical for the nonlinear coupling of these waves. The reason is mainly from the critical nonlinear convolution $K_3$ in Lemma \ref{A.5}, where an extra $(1+t)^{-\frac{1}{2}}$ than the H-wave above in the first part of the nonlinear convolution is usually from the conservation (or the divergence form of nonlinear terms). On the other hand, we know that the system (\ref{1.1}) is not conservative due to two damped terms in two momentum equations. All the same, the divergence form of the other nonlinear terms will still be used in treating the nonlinear coupling.
Of course, we should use different methods to these two kinds of nonlinear terms. Therefore, we consider the variables $\rho,\ m=\rho u,\ n,\ \omega =nv$, and the initial data for the system \eqref{1.1} is given as follows:
\begin{equation}\label{1.2}
(\rho,m,n,\omega)(x,t)|_{t=0}
 =(\rho_0,m_0,n_0,\omega_0)(x),\ \ x\in \mathbb{R}^3.
\end{equation}

The main novelty of the present paper is to develop the pointwise space-time description in $H^5$-framework, which is stated in the following theorem.
\begin{theorem}\label{l 1}Assume that $(\rho_0-\bar{\rho},m_0,n_0-\bar{n},\omega_0)\in H^5(\mathbb{R}^3)$ with  $\bar{\rho}>0$, $\bar{n}>0$ and $\varepsilon_0=:\!\!\|(\rho_0-\bar{\rho},m_0,n_0-\bar{n},\omega_0)\|_{H^5(\mathbb{R}^3)}$ small. Then there is a unique global classical solution $(\rho,m,n,\omega)$ of the Cauchy problem (\ref{1.1})-(\ref{1.2}). If further for $|\alpha|\leq 2$ and $|\tilde{\alpha}|\leq 1$,
\begin{equation}\label{1.3}
|D_x^\alpha(\rho_{0}-\bar{\rho},m_0)|\leq \varepsilon_0(1+|x|^2)^{-r},\ \
|D_x^{\tilde{\alpha}}(n_0-\bar{n},\omega_0)|\leq \varepsilon_0(1+|x|^2)^{-r},\ r>\frac{21}{10},
\end{equation}
then for the base sound speed $c:=\sqrt{\frac{\bar{n}P'(\bar{n})+\bar{\rho}}{\bar{n}+\bar{\rho}}}>0$, it holds for $|\tilde{\alpha}|\leq1$ that
\begin{equation}\label{1.4}
|D_x^{\tilde{\alpha}}(\rho-\bar{\rho},m,n-\bar{n},\omega)|\leq C(1+t)^{-\frac{4+|\tilde{\alpha}|}{2}}\Big(1+\frac{(|x|-ct)^2}{1+t}\Big)^{-\frac{3}{2}}\!
+C(1+t)^{-\frac{3+|\tilde{\alpha}|}{2}}\Big(1+\frac{|x|^2}{1+t}\Big)^{-\frac{3}{2}}.
\end{equation}
\end{theorem}

Note that the global existence and $L^2$-decay rate in $H^l$ with $l\geq3$ have been been given in \cite{Wugc}, our main contribution here is the pointwise space-time estimates as above. To the best of our knowledge, it is the first result in this direction for two-phase models.
\begin{remark}
The base sound speed $c:=\sqrt{\frac{\bar{n}P'(\bar{n})+\bar{\rho}}{\bar{n}+\bar{\rho}}}$ for the two-phase model (\ref{1.1}) is different from $c:=\sqrt{P'(\bar{n})}$ for the compressible isentropic Navier-Stokes system (NS). However, when $\bar{\rho}\ll\bar{n}$, it is almost the same as NS, which is naturally reasonable.
\end{remark}

\begin{remark}The pointwise estimates in Theorem \ref{l 1} for the Cauchy problem (\ref{1.1})-(\ref{1.2}) are similar to the compressible Navier-Stokes system in \cite{Du,LS,LW}, which also exhibit the generalized Huygens' principle. The byproduct, $L^p$-estimates also imply the dominant part of $(\rho-\bar{\rho},m,n-\bar{n},\omega)$ is the Huygens' wave when $p\leq2$ and the diffusion wave when $p\geq2$. In particular,
\bess
&&\|(\rho-\bar{\rho},m,n-\bar{n},\omega)(\cdot,t)\|_{L^p(\mathbb{R}^3)}\leq
\bigg\{\begin{array}{ll}
  C(1+t)^{-(2-\frac{5}{2p})}, & \ \ \ 1<p\leq2, \\[0.8mm]
 C(1+t)^{-\frac{3}{2}(1-\frac{1}{p})}, & \ \ \ 2\leq p\leq\infty.
\end{array}
\eess
This $L^p$-decay rate is a generalization of $L^2$-decay rate in \cite{Wugc}.
\end{remark}

Now, we introduce the main steps and difficulties in this paper. First of all, we use the Hodge decomposition for this big system to decompose the linear system for the densities and momentums into the compressible and incompressible parts, which is different from that for the densities and velocities in \cite{Wugc}. This will naturally impact the parameters in the Green's function, thus we should reconstruct the explicit representation of the Green's function in the Fourier space. Second, we shall give the spectral analysis by using low-frequency and high-frequency decomposition. Due to using the Hodge decomposition, we will meet the nonlocal Riesz operator with the symbol $\frac{\xi}{|\xi|}$, which formally brings the new singularity in the low frequency part. This singularity together with the wave operator and the heat operator in the low frequency part of the Green's function forces us to find suitable combinations such that we can derive desired the Huygens' wave and the diffusion wave of the Green's function in the physical space. On the other hand, due to the definition of the Green's function, there are the singularities in the high frequency part of the Green's function, and these singularities directly results in the regularity requirement of the initial data due to the quasi-linearity of the original system. To minimize the the regularity requirement, we have to carefully obtain the description of the singularity for each component in the high frequency part of the Green's function, and the details can be seen in Lemma \ref{l 3.4}. Last but not least, for the nonlinear coupling, non-conservative structure of (\ref{1.1}) from the damped mechanism also forces us to find some cancelations in the related columns of the Green's function to overcome the nonlinear coupling with these non-conservative terms (nonlinear terms without a divergence form). For the last difficulty, please see the details in (\ref{2.2}), Lemma \ref{l 3.1}, Lemma \ref{l 3.2} and Lemma \ref{A.5}.

\textbf{Notation.} We give some notations used in this paper. $C$ denotes a general positive constant which may vary in different estimates. We use $H^s(\mathbb{R}^n)=W^{s,2}(\mathbb{R}^n)$, where $W^{s,p}(\mathbb{R}^n)$ is the usual Sobolev space with its norm $\|f\|_{W^{s,p}(\mathbb{R}^n)}=\sum\limits_{k=0}^s\|\partial_x^kf\|_{L^p(\mathbb{R}^n)}$.

The remainder of the paper is organized as follows. In Section 2, we give the spectrum analysis of the linear system. Section 3 establishes the pointwise estimates for the Green's function. In Section 4, we deduce the pointwise estimates for the nonlinear system and prove Theorem \ref{l 1}. In the appendix, some useful inequalities are provided.

\section{Green's function}
\subsection{Linearization and Reformulation}
We first reformulate the system (\ref{1.1}). In what follows, we assume the steady state of the Cauchy problem \eqref{1.1}-\eqref{1.2} is $(\bar \rho, 0,\bar n, 0)$. For simplicity, we still use $(\rho,m,n,\omega)$ to denote the perturbation $(\rho-\bar{\rho},m,n-\bar{n},\omega)$ without confusion. Setting
\begin{equation*}
m=\rho u,\ \ w=\rho v,\ \ \alpha_{1}=P^{\prime}(\bar{n}),\ \ \alpha_{2}=\frac{\bar{\rho}}{\bar{n}},\ \ \bar{\mu}=\frac{\mu}{\bar{n}},\ \ \bar{\lambda}=\frac{\lambda}{\bar{n}},
\end{equation*}
then the system (\ref{1.1}) can be rewritten in the perturbation form as
\begin{equation}\label{2.1}
\left\{\begin{array}{l}
\rho_{t}+\operatorname{div} m=0, \\
m_{t}+\nabla \rho+m-\alpha_{2}w=F_{1}(\rho,m,n,\omega), \\
n_{t}+\operatorname{div}w=0, \\
w_{t}+\alpha_{1}\nabla n-\bar{\mu} \Delta w-(\bar{\mu}+\bar{\lambda}) \nabla \operatorname{div} w-m+\alpha_{2}w=F_{2}(\rho,m,n,\omega), \\
\left.(\rho, m, n, w)\right|_{t=0}=\left(\rho_{0}(x), m_{0}(x), n_{0}(x), w_{0}(x)\right),
\end{array}\right.
\end{equation}
where
\begin{equation}\label{2.2}
\begin{array}{rl}
F_{1}=&\D-\operatorname{div}\Big(\frac{m\otimes m}{\rho+\bar{\rho}}\Big)+\Big(\frac{\rho+\bar{\rho}}{n+\bar{n}}-\frac{\bar{\rho}}{\bar{n}}\Big)\omega,\\[2mm]
F_{2}=&\D-\operatorname{div}\Big(\frac{w\otimes m}{n+\bar{n}}\Big)-\bar{\mu}\Delta\frac{nw}{n+\bar{n}}-(\bar{\mu}+\bar{\lambda})\nabla\operatorname{div}\frac{nw}{n+\bar{n}}
-\nabla\Big(\rho(n+\bar{n})-\alpha_{1}n\Big)\\[2mm]
&\D-\Big(\frac{\rho+\bar{\rho}}{n+\bar{n}}-\frac{\bar{\rho}}{\bar{n}}\Big)\omega.
\end{array}
\end{equation}
Note that the last terms in $F_1$ and $F_2$ have not a divergence form, which bring us much more difficulties in deducing the generalized Huygens' principle.

Define $U=(\rho, m, n, w)^{T}$. In terms of the semigroup theory for evolutionary equation, we will study the following IVP for the linearized two-phase fluid system:
\begin{equation}\label{2.3}
\left\{\begin{array}{l}
U_{t}=\mathcal{A} U, \\
\left.U\right|_{t=0}=U_{0},
\end{array}\right.
\end{equation}
where the operator $\mathcal{A}$ is given by
\begin{equation*}
\mathcal{A}=\left(\begin{array}{cccc}
0 & -\operatorname{div} & 0 & 0 \\
-\nabla & -I_{3 \times 3} & 0 & \alpha_{2}I_{3 \times 3} \\
0 & 0 & 0 & -\operatorname{div} \\
0 & I_{3 \times 3} & -\alpha_{1} \nabla & \left(\bar{\mu} \Delta-\alpha_{2}\right) I_{3 \times 3}+(\bar{\mu}+\bar{\lambda}) \nabla \otimes \nabla
\end{array}\right).
\end{equation*}
Applying the Fourier transform to the system (\ref{2.3}), we have
\begin{equation}\label{2.4}
\left\{\begin{array}{l}
\hat{U}_{t}=\mathcal{A}(\xi) \hat{U}, \\
\left.\hat{U}\right|_{t=0}=\widehat{U}_{0},
\end{array}\right.
\end{equation}
where $\widehat{U}(\xi, t)=\mathcal{F}(U(x, t)), \xi=\left(\xi^{1}, \xi^{2}, \xi^{3}\right)^{T}$ and $\mathcal{A}(\xi)$ is defined by
\begin{equation}\label{2.5}
\mathcal{A}(\xi)=\left(\begin{array}{cccc}
0 & -i \xi^{T} & 0 & 0 \\
-i \xi & -I_{3 \times 3} & 0 & \alpha_{2}I_{3 \times 3} \\
0 & 0 & 0 & -i \xi^{T} \\
0 &  I_{3 \times 3} & -\alpha_{1} i \xi^{T} & -\left(\bar{\mu}|\xi|^{2}+\alpha_{2}\right) I_{3 \times 3}-(\bar{\mu}+\bar{\lambda}) \xi\xi^T
\end{array}\right).
\end{equation}
To facilitate narrative in the proof of the pointwise space-time estimates for the nonlinear problem in the last section, we also use the definition of the Green's function $G(x,t)$  with the following standard form as our previous works:
\begin{equation}\label{2.6}
\left\{\begin{array}{l}
G_{t}=\mathcal{A} G, \\
G|_{t=0}=\delta_0(x)I_8.
\end{array}\right.
\end{equation}

Besides, to give the representation of the Green's function in the Fourier space more easily, we use the Hodge decomposition. Let $\varphi=\Lambda^{-1} \operatorname{div} m$ and  $\psi=\Lambda^{-1} \operatorname{div} w$  be the ``compressible part" of the momenta $m$ and $w$, respectively, and denote  $\Phi=\Lambda^{-1} \operatorname{curl} m$ and  $\Psi=\Lambda^{-1} \operatorname{curl} w$  $( {\rm with} \ (\operatorname{curl} z)_{i}^{j}=   \left.\partial_{x_{j}} z^{i}-\partial_{x_{i}} z^{j}\right)$ by the ``incompressible part" of the momenta $m$ and $w$, respectively. Then, we the system (\ref{2.4}) becomes
\begin{equation}\label{2.7}
\left\{\begin{array}{l}
\rho_{t}+\Lambda \varphi=0, \\
\varphi_{t}-\Lambda \rho+(\varphi-\alpha_{2}\psi)=0, \\
n_{t}+\Lambda \psi=0,\\
\psi_{t}-\alpha_{1} \Lambda n+(2\bar{\mu}+\bar{\lambda}) \Lambda^{2} \psi-(\varphi-\alpha_{2}\psi)=0, \\
\left.(\rho, \varphi, n, \psi)\right|_{t=0}=\left(\rho_{0}(x), \Lambda^{-1} \operatorname{div} m_{0}(x), n_{0}(x), \Lambda^{-1} \operatorname{div} w_{0}(x)\right),
\end{array}\right.
\end{equation}
and
\begin{equation}\label{2.8}
\left\{\begin{array}{l}
\Phi_{t}+(\Phi-\alpha_{2}\Psi)=0, \\
\Psi_{t}+\mu \Lambda^{2} \Psi-(\Phi-\alpha_{2}\Psi)=0, \\
\left.(\Phi, \Psi)\right|_{t=0}=\left(\Lambda^{-1} \operatorname{curl} m_{0}(x), \Lambda^{-1} \operatorname{curl} w_{0}(x)\right).
\end{array}\right.
\end{equation}

\subsection{Spectral analysis for the compressible part}

We shall write the IVP  (\ref{2.7})  for  $\mathcal{U}=(\rho, \varphi, n, \psi)^{t} $ as
\begin{equation}\label{2.9}
\left\{\begin{array}{l}
\mathcal{U}_{t}=\mathcal{A}_{1} \mathcal{U}, \\
\left.\mathcal{U}\right|_{t=0}=\mathcal{U}_{0},
\end{array}\right.
\end{equation}
where the operator  $\mathcal{A}_{1}$  is given by
\begin{equation*}
\mathcal{A}_{1}=\left(\begin{array}{cccc}
0 & -\Lambda & 0 & 0 \\
\Lambda & -1 & 0 & \alpha_{2} \\
0 & 0 & 0 & - \Lambda \\
0 & 1 & \alpha_{1} \Lambda & -\nu \Lambda^{2}-\alpha_{2}
\end{array}\right),\ {\rm with}\ \ \nu=2 \bar{\mu}+\bar{\lambda}.
\end{equation*}
Taking the Fourier transform to the system (\ref{2.9}), we have
\begin{equation*}
\left\{\begin{array}{l}
\hat{\mathcal{U}}_{t}=\mathcal{A}_{1}(\xi) \hat{\mathcal{U}}, \\
\left.\hat{\mathcal{U}}\right|_{t=0}=\hat{\mathcal{U}}_{0},
\end{array}\right.
\end{equation*}
where  $\hat{\mathcal{U}}(\xi, t)=\mathcal{F}(\mathcal{U}(x, t)) $ and  $\mathcal{A}_{1}(\xi)$  is defined by
\begin{equation}
\mathcal{A}_{1}(\xi)=\left(\begin{array}{cccc}
0 & -|\xi| & 0 & 0 \\
|\xi| & -1 & 0 & \alpha_{2} \\
0 & 0 & 0 & -|\xi| \\
0 & 1 & \alpha_{1}|\xi| & -\nu|\xi|^{2}-\alpha_{2}
\end{array}\right).
\end{equation}
Its eigenvalues satisfy
\begin{equation}\label{2.10}
\begin{array}{rl}
&{\rm det}\left(r \mathrm{I}-\mathcal{A}_{1}(\xi)\right) \\
=& r^{4}+\left[\nu|\xi|^{2}+\alpha_{2}+1\right] r^{3}+\left( \nu+\alpha_{1} +1\right)|\xi|^{2} r^{2} \\
&+\left[\nu|\xi|^{4}+\left(\alpha_{1}+\alpha_{2}\right)|\xi|^{2}\right] r+\alpha_{1} |\xi|^{4}
=0.
\end{array}
\end{equation}
Therefore the matrix  $\mathcal{A}_{1}(\xi)$  has four different eigenvalues: $r_{1}(|\xi|), r_{2}(|\xi|), r_{3}(|\xi|), r_{4}(|\xi|)$.
Thus, the semigroup $\mathrm{e}^{t \mathcal{A}_{1}}$ can be decomposed into
\begin{equation}\label{2.11}
\mathrm{e}^{t \mathcal{A}_{1}(\xi)}=\sum_{i=1}^{4} \mathrm{e}^{r_{i} t} P_{i}(\xi),
\end{equation}
and the projector  $P_{i}(\xi)$  is
\begin{equation}\label{2.12}
P_{i}(\xi)=\prod_{j \neq i} \frac{\mathcal{A}_{1}(\xi)-r_{j} I}{r_{i}-r_{j}}, \quad i, j=1,2,3,4.
\end{equation}
Therefore, we can show the solution of IVP (\ref{2.7}) as
\begin{equation}\label{2.13}
\hat{\mathcal{U}}(\xi, t)=\mathrm{e}^{t \mathcal{A}_{1}(\xi)} \hat{\mathcal{U}}_{0}(\xi)=\bigg(\sum_{i=1}^{4} \mathrm{e}^{r_{i} t} P_{i}(\xi)\bigg) \widehat{\mathcal{U}}_{0}(\xi).
\end{equation}

Here and below, we use the superscript $``l"$ to denote the low frequency part, and use the superscript $``h"$ means the high frequency part.

\vspace{4mm}
\textbf{\textit{Low frequency part.}}\vspace{3mm}

By a direct computation, we have the following for the spectral in the low frequency part:
\begin{lemma}\label{l 2.1} There exists a positive constant  $\eta_{1} \ll 1$  such that, for  $|\xi| \leq \eta_{1}$, the spectral has the following Taylor series expansion:
\begin{equation}\label{2.16}
\left\{\begin{aligned}
r_{1}=&-\alpha_{2}-1+\frac{-\alpha_{2}\left(\alpha_{2}+1\right)\nu+\alpha_{1} \alpha_{2} +1}{\left(\alpha_{2}+1\right)^{2}}|\xi|^{2}+O\left(|\xi|^{4}\right), \\
r_{2}=&-\frac{\alpha_{1}}{\alpha_{1}+\alpha_{2}}|\xi|^{2}+O\left(|\xi|^{4}\right), \\
r_{3}=&-\frac{\nu\left(\alpha_{1}+\alpha_{2}\right)\left(\alpha_{2}+1\right)+\alpha_{2}\left(\alpha_{1}-1\right)^{2}}{2\left(\alpha_{1} +\alpha_{2}\right)\left(\alpha_{2}+1\right)^{2}}|\xi|^{2}+O\left(|\xi|^{4}\right)\\
&+ i\left[\sqrt{\frac{\alpha_{1}+\alpha_{2}}{\alpha_{2}+1}}|\xi|+O\left(|\xi|^{3}\right)\right], \\
r_{4}=&-\frac{\nu\left(\alpha_{1}+\alpha_{2}\right)\left(\alpha_{2}+1\right)+\alpha_{2}\left(\alpha_{1} -1\right)^{2}}{2\left(\alpha_{1} +\alpha_{2}\right)\left(\alpha_{2}+1\right)^{2}}|\xi|^{2}+O\left(|\xi|^{4}\right) \\
&- i\left[\sqrt{\frac{\alpha_{1} +\alpha_{2}}{\alpha_{2}+1}}|\xi|+O\left(|\xi|^{3}\right)\right].
\end{aligned}\right.
\end{equation}
\end{lemma}

 Therefore, for $|\xi|\ll1$ it holds that
\begin{equation}
\begin{aligned}
P_{1}^l(\xi)=&\frac{\mathcal{A}_1-r_{2}I}{r_{1}-r_{2}} \frac{\mathcal{A}_{1}-r_{3}I}{r_{1}-r_{3}} \frac{\mathcal{A}_{1}-r_{4}I}{r_{1}-r_{4}}
=&-\frac{1}{\left(\alpha_{2}+1\right)}\left(\begin{array}{cccc}
0 & 0 & 0 & 0 \\
0 & -1 & 0 & \alpha_{2} \\
0 & 0 & 0 & 0 \\
0 & 1& 0 & -\alpha_{2}
\end{array}\right)+\mathcal{O}(|\xi|),
\end{aligned}
\end{equation}
\begin{equation}
\begin{aligned}
P_{2}^l(\xi)=&\frac{\mathcal{A}_{1}-r_{1}I}{r_{2}-r_{1}} \frac{\mathcal{A}_{1}-r_{3}I}{r_{2}-r_{3}} \frac{\mathcal{A}_{1}-r_{4}I}{r_{2}-r_{4}}
=&\frac{1}{\left(\alpha_{1}+\alpha_{2}\right)}\left(\begin{array}{cccc}
\alpha_{1}  & 0 & -\alpha_{1}\alpha_{2} & 0 \\
0 & 0 & 0 & 0 \\
-1  & 0 & \alpha_{2} & 0 \\
0 & 0 & 0 & 0
\end{array}\right)+\mathcal{O}(|\xi|),
\end{aligned}
\end{equation}
\begin{equation}
\begin{aligned}
&P_{3}^l(\xi)=\frac{\mathcal{A}_{1}-r_{1}I}{r_{3}-r_{1}} \frac{\mathcal{A}_{1}-r_{2}I}{r_{3}-r_{2}} \frac{\mathcal{A}_{1}-r_{4}I}{r_{3}-r_{4}}\\
&=\frac{-1}{2\left(\alpha_{1} +\alpha_{2}\right)}\left(\begin{array}{cccc}
-\alpha_{2} & -i \alpha_{2} \sqrt{\frac{\alpha_{1} +\alpha_{2}}{\alpha_{2}+1}} & -\alpha_{1}\alpha_{2} & -i\alpha_{2} \sqrt{\frac{\alpha_{1} +\alpha_{2}}{\alpha_{2}+1}} \\
i \alpha_{2} \sqrt{\frac{\alpha_{1}+\alpha_{2}}{\alpha_{2}+1}} & -\frac{\alpha_{2}\left(\alpha_{1}+\alpha_{2}\right)}{\alpha_{2}+1} & i \alpha_{1}\alpha_{2} \sqrt{\frac{\alpha_{1} +\alpha_{2}}{\alpha_{2}+1}} & -\frac{\alpha_{2}\left(\alpha_{1} +\alpha_{2}\right)}{\alpha_{2}+1} \\
-1 & -i \sqrt{\frac{\alpha_{1}+\alpha_{2}}{\alpha_{2}+1}} & -\alpha_{1} & -i \sqrt{\frac{\alpha_{1} +\alpha_{2}}{\alpha_{2}+1}} \\
i \sqrt{\frac{\alpha_{1}+\alpha_{2}}{\alpha_{2}+1}} & -\frac{\alpha_{1} +\alpha_{2}}{\alpha_{2}+1} & i \alpha_{1} \sqrt{\frac{\alpha_{1} +\alpha_{2}}{\alpha_{2}+1}} & -\frac{\alpha_{1} +\alpha_{2}}{\alpha_{2}+1}
\end{array}\right)+\mathcal{O}(|\xi|),
\end{aligned}
\end{equation}
\begin{equation}
\begin{aligned}
&P_{4}^l(\xi)=\frac{\mathcal{A}_{1}-r_{1}I}{r_{4}-r_{1}} \frac{\mathcal{A}_{1}-r_{2}I}{r_{4}-r_{2}} \frac{\mathcal{A}_{1}-r_{3}I}{r_{4}-r_{3}}\\
=&\frac{-1}{2\left(\alpha_{1} +\alpha_{2}\right)}\left(\begin{array}{cccc}
-\alpha_{2} & i \alpha_{2} \sqrt{\frac{\alpha_{1} +\alpha_{2}}{\alpha_{2}+1}} & -\alpha_{1}\alpha_{2} & i\alpha_{2} \sqrt{\frac{\alpha_{1} +\alpha_{2}}{\alpha_{2}+1}} \\
-i\alpha_{2} \sqrt{\frac{\alpha_{1}+\alpha_{2}}{\alpha_{2}+1}} & -\frac{\alpha_{2}\left(\alpha_{1}+\alpha_{2}\right)}{\alpha_{2}+1} & -i \alpha_{1}\alpha_{2} \sqrt{\frac{\alpha_{1} +\alpha_{2}}{\alpha_{2}+1}} & -\frac{\alpha_{2}\left(\alpha_{1} +\alpha_{2}\right)}{\alpha_{2}+1} \\
-1 & i \sqrt{\frac{\alpha_{1}+\alpha_{2}}{\alpha_{2}+1}} & -\alpha_{1} & i \sqrt{\frac{\alpha_{1} +\alpha_{2}}{\alpha_{2}+1}} \\
-i\sqrt{\frac{\alpha_{1}+\alpha_{2}}{\alpha_{2}+1}} & -\frac{\alpha_{1} +\alpha_{2}}{\alpha_{2}+1} & -i \alpha_{1} \sqrt{\frac{\alpha_{1} +\alpha_{2}}{\alpha_{2}+1}} & -\frac{\alpha_{1} +\alpha_{2}}{\alpha_{2}+1}
\end{array}\right)+\mathcal{O}(|\xi|).
\end{aligned}
\end{equation}
Here and below, we use $``\cdots"$ to denote the rest terms, and the rest terms don't impact the estimates. Substituting the above estimates in (\ref{2.13}), one has:
\begin{lemma}\label{l 2.2} There exists a positive constant  $\eta_{1} \ll 1$  such that, for  $|\xi| \leq \eta_{1}$, we can induce
\begin{equation}\label{2.21}
\begin{aligned}
\hat{\rho}^{l}=&\frac{\mathrm{e}^{r_{2} t}}{2(\alpha_{1}+\alpha_{2})}\big(2\alpha_{1} \hat{\rho}_{0}^{l}-2\alpha_{1}\alpha_{2} \hat{n}_{0}^{l} \big)\\
&-\frac{\alpha_{2}\mathrm{e}^{r_{3} t}}{2(\alpha_{1}+\alpha_{2})}\bigg(- \hat{\rho}_{0}^{l}-i  \sqrt{\frac{\alpha_{1} +\alpha_{2}}{\alpha_{2}+1}}  \hat{\varphi}_{0}^{l}-\alpha_{1} \hat{n}_{0}^{l}-i\sqrt{\frac{\alpha_{1} +\alpha_{2}}{\alpha_{2}+1}} \hat{\psi}_{0}^{l} \bigg)\\
&-\frac{\alpha_{2}\mathrm{e}^{r_{4} t}}{2(\alpha_{1}+\alpha_{2})}\bigg(-\hat{\rho}_{0}^{l}+i \sqrt{\frac{\alpha_{1} +\alpha_{2}}{\alpha_{2}+1}}  \hat{\varphi}_{0}^{l}-\alpha_{1} \hat{n}_{0}^{l}+i \sqrt{\frac{\alpha_{1} +\alpha_{2}}{\alpha_{2}+1}} \hat{\psi}_{0}^{l} \bigg)+\cdots,
\end{aligned}
\end{equation}
\begin{equation}\label{2.22}
\begin{aligned}
\hat{n}^{l}=&\frac{\mathrm{e}^{r_{2} t}}{2(\alpha_{1}+\alpha_{2})}\big(-2 \hat{\rho}_{0}^{l}+\alpha_{2} \hat{n}_{0}^{l}\big)\\
&-\frac{\mathrm{e}^{r_{3} t}}{2\big(\alpha_{1}+\alpha_{2}\big)}\big(-\hat{\rho}_{0}^{l}-i \sqrt{\frac{\alpha_{1} +\alpha_{2}}{\alpha_{2}+1}}  \hat{\varphi}_{0}^{l}-\alpha_{1}\hat{n}_{0}^{l}-i\sqrt{\frac{\alpha_{1} +\alpha_{2}}{\alpha_{2}+1}} \hat{\psi}_{0}^{l} \big)\\
&-\frac{\mathrm{e}^{r_{4} t}}{2(\alpha_{1}+\alpha_{2})}\big(-\alpha_{2} \hat{\rho}_{0}^{l}+i \sqrt{\frac{\alpha_{1} +\alpha_{2}}{\alpha_{2}+1}}  \hat{\varphi}_{0}^{l}-\alpha_{1} \hat{n}_{0}^{l}+i \sqrt{\frac{\alpha_{1} +\alpha_{2}}{\alpha_{2}+1}} \hat{\psi}_{0}^{l} \big)+\cdots,
\end{aligned}
\end{equation}
\begin{equation}\label{2.23}
\begin{aligned}
\hat{\varphi}^{l}=&-\frac{\mathrm{e}^{r_{1} t}}{\alpha_{2}+1}\big(- \hat{\varphi}_{0}^{l}+ \alpha_{2}\hat{\psi}_{0}^{l}\big)\\
&-\frac{\alpha_{2}\mathrm{e}^{r_{3} t}}{2(\alpha_{1}+\alpha_{2})}\big(i \sqrt{\frac{\alpha_{1} +\alpha_{2}}{\alpha_{2}+1}} \hat{\rho}_{0}^{l}
- \frac{(\alpha_{1} +\alpha_{2})}{\alpha_{2}+1} \hat{\varphi}_{0}^{l}+i\alpha_{1}\sqrt{\frac{\alpha_{1} +\alpha_{2}}{\alpha_{2}+1}} \hat{n}_{0}^{l}
- \frac{ \alpha_{1} +\alpha_{2}}{\alpha_{2}+1}  \hat{\psi}_{0}^{l}\big)\\
&+\frac{\alpha_{2}\mathrm{e}^{r_{4} t}}{2(\alpha_{1}+\alpha_{2})}(i \sqrt{\frac{\alpha_{1} +\alpha_{2}}{\alpha_{2}+1}} \hat{\rho}_{0}^{l}
+\frac{ (\alpha_{1} +\alpha_{2})}{\alpha_{2}+1} \hat{\varphi}_{0}^{l}+i\alpha_{1}\sqrt{\frac{\alpha_{1} +\alpha_{2}}{\alpha_{2}+1}} \hat{n}_{0}^{l}+\frac{ \alpha_{1} +\alpha_{2}}{\alpha_{2}+1}  \hat{\psi}_{0}^{l}\big)+\cdots,\\
\sim & \frac{1}{2(\alpha_{1}+\alpha_{2})}\big(-i \alpha_{2} \sqrt{\frac{\alpha_{1}+\alpha_{2}}{\alpha_{2}+1}}\mathrm{e}^{r_{3}t}+i \alpha_{2} \sqrt{\frac{\alpha_{1} +\alpha_{2}}{\alpha_{2}+1}}\mathrm{e}^{r_{4} t} \big)\hat{\rho}_{0}^{l}\\
&+\frac{1}{2(\alpha_{2}+1)}\big(2\alpha_2\mathrm{e}^{r_{1} t}+\alpha_{2} \mathrm{e}^{r_{3}t}+\alpha_{2} \mathrm{e}^{r_{4}t}\big)i\frac{\xi \hat{m}_{0}^{l}}{|\xi |}\\
&+\frac{1}{2(\alpha_{1}+\alpha_{2})}\big(-i \alpha_{1} \alpha_{2}\sqrt{\frac{\alpha_{1}+\alpha_{2}}{\alpha_{2}+1}}\mathrm{e}^{r_{3}t}+i \alpha_{1}\alpha_{2} \sqrt{\frac{\alpha_{1} +\alpha_{2}}{\alpha_{2}+1}}\mathrm{e}^{r_{4} t}\big)\hat{n}_{0}^{l}\\
&+\frac{1}{2(\alpha_{2}+1)}\big(-2\alpha_{2}\mathrm{e}^{r_{1} t}+\alpha_{2}\mathrm{e}^{r_{3}t}+ \alpha_{2}\mathrm{e}^{r_{4}t}\big)i\frac{\xi \hat{\omega}_{0}^{l}}{|\xi |}+\cdots,
\end{aligned}
\end{equation}
\begin{equation}\label{2.24}
\begin{aligned}
\hat{\psi}^{l}=&-\frac{\mathrm{e}^{r_{1} t}}{\alpha_{2}+1}\big( \hat{\varphi}_{0}^{l}-\alpha_{2} \hat{\psi}_{0}^{l}\big)\\
&-\frac{\mathrm{e}^{r_{3} t}}{2(\alpha_{1}+\alpha_{2})}(i  \sqrt{\frac{\alpha_{1} +\alpha_{2}}{\alpha_{2}+1}} \hat{\rho}_{0}^{l}
- \frac{ (\alpha_{1} +\alpha_{2})}{\alpha_{2}+1} \hat{\varphi}_{0}^{l}+i\alpha_{1}\sqrt{\frac{\alpha_{1}+\alpha_{2}}{\alpha_{2}+1}} \hat{n}_{0}^{l}
- \frac{ \alpha_{1} +\alpha_{2}}{\alpha_{2}+1}  \hat{\psi}_{0}^{l})\\
&+\frac{\mathrm{e}^{r_{4} t}}{2(\alpha_{1}+\alpha_{2})}(i  \sqrt{\frac{\alpha_{1} +\alpha_{2}}{\alpha_{2}+1}} \hat{\rho}_{0}^{l}
+\frac{ (\alpha_{1} +\alpha_{2})}{\alpha_{2}+1} \hat{\varphi}_{0}^{l}+i\alpha_{1}\sqrt{\frac{\alpha_{1} +\alpha_{2}}{\alpha_{2}+1}} \hat{n}_{0}^{l}
+\frac{ \alpha_{1} +\alpha_{2}}{\alpha_{2}+1}  \hat{\psi}_{0}^{l})+\cdots,\\
\sim &\frac{1}{2(\alpha_{1}+\alpha_{2})}\big(-i  \sqrt{\frac{\alpha_{1}+\alpha_{2}}{\alpha_{2}+1}}\mathrm{e}^{r_{3}t}+i  \sqrt{\frac{\alpha_{1} +\alpha_{2}}{\alpha_{2}+1}}\mathrm{e}^{r_{4} t}\big)\hat{\rho}_{0}^{l}\\
&+\frac{1}{2(\alpha_{2}+1)}(-2\mathrm{e}^{r_{1} t}+ \mathrm{e}^{r_{3}t}+ \mathrm{e}^{r_{4}t})i\frac{\xi \hat{m}_{0}^{l}}{|\xi |}\\
&+\frac{1}{2(\alpha_{1}+\alpha_{2})}\big(-i \alpha_{1} \sqrt{\frac{\alpha_{1}+\alpha_{2}}{\alpha_{2}+1}}\mathrm{e}^{r_{3}t}+i \alpha_{1} \sqrt{\frac{\alpha_{1} \bar{n}+\alpha_{2}}{\alpha_{2}+1}}\mathrm{e}^{r_{4} t}\big)\hat{n}_{0}^{l}\\
&+\frac{1}{2(\alpha_{2}+1)}(2\alpha_{2}\mathrm{e}^{r_{1} t}+\mathrm{e}^{r_{3}t}+ \mathrm{e}^{r_{4}t})i\frac{\xi \hat{\omega}_{0}^{l}}{|\xi |}+\cdots.
\end{aligned}
\end{equation}
\end{lemma}

\textbf{\textit{High frequency part.}}\vspace{3mm}

Then we consider the high frequency part. After a direct computation, one has
\begin{lemma}\label{l 2.3}
 There exists a positive constant  $\eta_{2} \gg 1$  such that, for  $|\xi| \gg \eta_{2}$, the spectral has the following Taylor series expansion:
\begin{equation}\label{2.24(1)}
\left\{\begin{aligned}
r_{1}=&-\frac{\alpha_{1}}{\nu}+O\left(|\xi|^{-2}\right),\\
r_{2}=&-\nu|\xi|^{2}+\frac{\alpha_{1}}{\nu}-\alpha_{2}+O\left(|\xi|^{-2}\right), \\
r_{3}=&-\frac{1}{2}+O\left(|\xi|^{-2}\right)+i|\xi|+i O\left(|\xi|^{-1}\right), \\
r_{4}=& -\frac{1}{2}+O\left(|\xi|^{-2}\right)-i|\xi|+i O\left(|\xi|^{-1}\right). \\
\end{aligned}\right.
\end{equation}
\end{lemma}
Substituting the estimates in Lemma \ref{l 2.3} into $P_i$ with $i=1,2,3,4$, one has:
\begin{lemma}\label{l 2.4}  There exists a positive constant  $\eta_{1} \gg 1$  such that, for  $|\xi| \gg \eta_{1}$, we can express  $P_{i}\ (1 \leq i \leq 4)$ as follows:
\begin{equation}\label{2.25}
\begin{aligned}
&P_{1}^h(\xi)=\frac{\mathcal{A}_{1}-r_{2}I}{r_{1}-r_{2}} \frac{\mathcal{A}_{1}-r_{3}I}{r_{1}-r_{3}} \frac{\mathcal{A}_{1}-r_{4}I}{r_{1}-r_{4}}\\
=&\frac{1}{\nu|\xi|^{4}}\left(\!\!\begin{array}{cccc}
\left(\frac{\alpha_{1}}{\nu}-\alpha_{2}\right)|\xi|^{2} &-\left(\frac{\alpha_{1}}{\nu}\right)^{2}|\xi|  & -\alpha_{1}\alpha_{2}|\xi|^{2} & \frac{\alpha_{1}\alpha_{2}}{\nu}|\xi| \\
\left(\frac{\alpha_{1}}{\nu}\right)^{2}|\xi|  & \left(\frac{\alpha_{1}}{\nu}-\alpha_{2}\right)|\xi|^{2} & -\frac{\alpha_{1}^{2}\alpha_{2}}{\nu}|\xi| & -\alpha_{2}\left((1+\alpha_{2})+\frac{\alpha_{1}\alpha_{2}}{\nu}+\left(\frac{\alpha_{1}}{\nu}\right)^{2}\right) \\
-|\xi|^{2} & \frac{\alpha_{1}}{\nu}|\xi| & \nu|\xi|^{4}  & -|\xi|^{3} \\
-\frac{\alpha_{1}}{\nu}|\xi|&\frac{\alpha_{1}}{\nu}\left(\alpha_{2}+\frac{\alpha_{1}}{\nu}\right)& \alpha_{1}|\xi|^{3} & -\alpha_{2}(1+\nu+\alpha_{1})|\xi|^{2}
\end{array}\!\!\right)+\cdots,
\end{aligned}
\end{equation}
\begin{equation*}
\begin{aligned}
&P_{2}^h(\xi)=\frac{\mathcal{A}_{1}-r_{1}I}{r_{2}-r_{1}} \frac{\mathcal{A}_{1}-r_{3}I}{r_{2}-r_{3}} \frac{\mathcal{A}_{1}-r_{4}I}{r_{2}-r_{4}}\\
=&\frac{-1}{\nu^{3}|\xi|^{6}}\left(\begin{array}{cccc}
0 &-\frac{\alpha_{1}}{\nu}\left(\frac{\alpha_{1}}{\nu}-\alpha_{2}\right)|\xi| & -\alpha_{1}\alpha_{2}|\xi|^{2} & \alpha_{2}\nu|\xi|^{3} \\
-\frac{\alpha_{1}}{\nu}\left(\frac{\alpha_{1}}{\nu}-\alpha_{2}\right)|\xi|&-\alpha_{2}\nu|\xi|^{2} & -\alpha_{1}\alpha_{2}\nu|\xi|^{3} & \alpha_{2}\nu^{2}|\xi|^{4} \\
-|\xi|^{2} & \nu|\xi|^{3} & \alpha_{1}\nu|\xi|^{4}  & -\nu^{2}|\xi|^{5} \\
-\nu|\xi|^{3}&\nu^{2}|\xi|^{4}& \alpha_{1}\nu^{2}|\xi|^{5} & -\nu^{3}|\xi|^{6}
\end{array}\right)+\cdots,
\end{aligned}
\end{equation*}
\begin{equation*}
\begin{aligned}
&P_{3}^h(\xi)=\frac{\mathcal{A}_{1}-r_{1}I}{r_{3}-r_{1}} \frac{\mathcal{A}_{1}-r_{2}I}{r_{3}-r_{2}} \frac{\mathcal{A}_{1}-r_{4}I}{r_{3}-r_{4}}\\
=&\frac{1}{4\nu^{2}|\xi|^{8}}\left(\begin{array}{cccc}
 2\nu^{2}|\xi|^{8}& i 2\nu^{2}|\xi|^{8} &  2\alpha_{1}\alpha_{2}\nu|\xi|^{6} & i2\alpha_{2}\nu|\xi|^{6} \\
-i 2\nu^{2}|\xi|^{8}& 2\nu^{2}|\xi|^{8} & -i2\alpha_{1}\alpha_{2}\nu|\xi|^{6} & 2\alpha_{2}\nu|\xi|^{6} \\
2\nu|\xi|^{6} &i 2\nu|\xi|^{6} & 0  & -i 2\alpha_{2}|\xi|^{6} \\
 -2\nu|\xi|^{6}& i2\nu|\xi|^{5}& -2\alpha_{1}\alpha_{2}\nu|\xi|^{7}& -\frac{1}{2}\nu^{2}|\xi|^{6}
\end{array}\right)+\cdots,
\end{aligned}
\end{equation*}
and
\begin{equation*}
\begin{aligned}
&P_{4}^h(\xi)=\frac{\mathcal{A}_{1}-r_{1}I}{r_{4}-r_{1}} \frac{\mathcal{A}_{1}-r_{2}I}{r_{4}-r_{2}} \frac{\mathcal{A}_{1}-r_{3}I}{r_{1}-r_{3}}\\
=&\frac{1}{4\nu^{2}|\xi|^{8}}\left(\begin{array}{cccc}
 2\nu^{2}|\xi|^{8}& -i 2\nu^{2}|\xi|^{8} &  2\alpha_{1}\alpha_{2}\nu|\xi|^{6} & -i2\alpha_{2}\nu|\xi|^{6} \\
 i 2\nu^{2}|\xi|^{8}& 2\nu^{2}|\xi|^{8} & i2\alpha_{1}\alpha_{2}\nu|\xi|^{6} & 2\alpha_{2}\nu|\xi|^{6} \\
2\nu|\xi|^{6} &-i 2\nu|\xi|^{6} & 0  & i 2\alpha_{2}|\xi|^{6} \\
 -2\nu|\xi|^{6}&-i2\nu|\xi|^{5}& -2\alpha_{1}\alpha_{2}\nu|\xi|^{7}& -\frac{1}{2}\nu^{2}|\xi|^{6}
\end{array}\right)+\cdots.
\end{aligned}
\end{equation*}
\end{lemma}
Additionally, we have the following asymptotic expansion of the unknowns in the Fourier space when $|\xi|\gg1$:
\begin{lemma}\label{l 2.5} There exists a positive constant  $\eta_{2} \gg 1$  such that, for  $|\xi| \gg \eta_{2}$, we can induce
\begin{equation}\label{2.27}
\begin{aligned}
\hat{\rho}^{h}
\sim &\Big( \frac{1}{\nu^{2}|\xi|^{2}}(\frac{\alpha_{1}}{\nu}-\alpha_{2})\mathrm{e}^{r_{1}t}+\frac{1}{2}\mathrm{e}^{r_{3} t} +\frac{1}{2}\mathrm{e}^{r_{4} t} \Big)\hat{\rho}_{0}^{h}\\
&+\Big(-\frac{\alpha_1^{2}}{\nu^{3}|\xi|^{3}}\mathrm{e}^{r_{1}t}+\frac{\alpha_1}{\nu^{4}|\xi|^{5}}\Big(\frac{\alpha_{1}}{\nu}-\alpha_{2}\Big)
\mathrm{e}^{r_{2}t} +i\frac{1}{2}\mathrm{e}^{r_{3} t} -i\frac{1}{2}\mathrm{e}^{r_{4} t}\Big)i\frac{\xi \hat{m}_{0}^{h}}{|\xi |}\\
&+\Big(-\frac{\alpha_1\alpha_{2}}{\nu|\xi|^{2}}\mathrm{e}^{r_{1}t}+\frac{\alpha_1\alpha_{2}}{\nu^{3}|\xi|^{4}}\mathrm{e}^{r_{2}t} +\frac{\alpha_1\alpha_{2}}{2\nu|\xi|^{2}}\mathrm{e}^{r_{3}t}+\frac{\alpha_1\alpha_{2}}{2\nu|\xi|^{2}}\mathrm{e}^{r_{4}t} \Big)\hat{n}_{0}^{h}\\
&+\Big(\frac{\alpha_1\alpha_{2}}{\nu^{2}|\xi|^{3}}\mathrm{e}^{r_{1}t}-\frac{\alpha_{2}}{\nu^{2}|\xi|^{3}}\mathrm{e}^{r_{2}t} +i\frac{\alpha_{2}}{2\nu|\xi|^{2}}\mathrm{e}^{r_{3}t}-i\frac{\alpha_{2}}{2\nu|\xi|^{2}}\mathrm{e}^{r_{4}t} \Big)i\frac{\xi\hat{\omega}_{0}^{h}}{|\xi |}+\cdots,\ \ \ \ \
\end{aligned}
\end{equation}
\begin{equation}\label{2.28}
\begin{aligned}
\hat{n}^{h}
\sim&\Big(-\frac{1}{\nu|\xi|^{2}}\mathrm{e}^{r_{1}t}+\frac{1}{\nu^{3}|\xi|^{4}}\mathrm{e}^{r_{2}t}+\frac{1}{2\nu|\xi|^{2}}\mathrm{e}^{r_{3}t}
+\frac{1}{2\nu|\xi|^{2}}\mathrm{e}^{r_{4}t}\Big)\hat{\rho}_{0}^{h}\\
&+\Big(\frac{\alpha_{1}}{\nu^{2}|\xi|^{3}}\mathrm{e}^{r_{1}t}-\frac{1}{\nu^{2}|\xi|^{3}}\mathrm{e}^{r_{2}t}
+i\frac{1}{2\nu|\xi|^{2}}\mathrm{e}^{r_{3}t}-i\frac{1}{2\nu|\xi|^{2}}\mathrm{e}^{r_{4}t}\Big)i\frac{\xi \hat{m}_{0}^{h}}{|\xi |} \\
&+\Big(\mathrm{e}^{r_{1}t}-\frac{\alpha_{1}}{\nu^{2}|\xi|^{2}}\mathrm{e}^{r_{2}t}\Big)\hat{n}_{0}^{h}
+\Big(-\frac{\mathrm{e}^{r_{1}t}-\mathrm{e}^{r_{2}t}}{\nu|\xi|}
-i\frac{\alpha_{2}(\mathrm{e}^{r_{3}t}-\mathrm{e}^{r_{4}t})}{2\nu^{2}|\xi|^{2}}\Big)
i\frac{\xi \hat{\omega}_{0}^{h}}{|\xi |}+\cdots,\ \ \
\end{aligned}
\end{equation}
\begin{equation}\label{2.29}
\begin{aligned}
\hat{\varphi}^{h}
\sim & \Big(\frac{\alpha_1^{2}}{\nu^{3}|\xi|^{3}}\mathrm{e}^{r_{1}t}-\frac{\alpha_1}{\nu^{4}|\xi|^{5}}\Big(\frac{\alpha_{1}}{\nu}-\alpha_{2}\Big)
\mathrm{e}^{r_{2}t} -i\frac{1}{2}\mathrm{e}^{r_{3} t} +i\frac{1}{2}\mathrm{e}^{r_{4} t}\Big)\hat{\rho}_{0}^{h}\\
&+\Big(\frac{1}{\nu^{2}|\xi|^{2}}(\frac{\alpha_{1}}{\nu}-\alpha_{2})\mathrm{e}^{r_{1}t}+\frac{\alpha_{2}}{\nu^{2}|\xi|^{4}}\mathrm{e}^{r_{2}t} +\frac{1}{2}\mathrm{e}^{r_{3} t} +\frac{1}{2}\mathrm{e}^{r_{4} t}\Big)i\frac{\xi \hat{m}_{0}^{h}}{|\xi |}\\
&+\Big(-\frac{\alpha_1^{2}\alpha_{2}}{\nu^{2}|\xi|^{3}}\mathrm{e}^{r_{1}t}+\frac{\alpha_{1}\alpha_{2}}{\nu^{2}|\xi|^{3}}\mathrm{e}^{r_{2}t} -i\frac{\alpha_{1}\alpha_{2}}{2\nu|\xi|^{2}}\mathrm{e}^{r_{3}t}+i\frac{\alpha_{1}\alpha_{2}}{2\nu|\xi|^{2}}\mathrm{e}^{r_{4}t} \Big)\hat{n}_{0}^{h}\\
&+\bigg(\Big(\frac{\alpha_{1}\alpha_{2}}{\nu}+(\frac{\alpha_{1}}{\nu})^{2}-1-\alpha_{2}\Big)\frac{\alpha_{2}\mathrm{e}^{r_{1}t}}{\nu|\xi|^{4}}
-\frac{\alpha_{2}\mathrm{e}^{r_{2}t} }{\nu|\xi|^{2}} +\frac{\alpha_{2}\mathrm{e}^{r_{3}t}}{2\nu|\xi|^{2}}+\frac{\mathrm{e}^{r_{4}t}}{2\nu|\xi|^{2}}\bigg)i\frac{\xi\hat{\omega}_{0}^{h}}{|\xi |}+\cdots,
\end{aligned}
\end{equation}
\begin{equation}\label{2.30}
\begin{aligned}
\hat{\psi}^{h}
\sim&\Big(-\frac{\alpha_{1}}{\nu^{2}|\xi|^{3}}\mathrm{e}^{r_{1}t}+\frac{1}{\nu^{2}|\xi|^{3}}\mathrm{e}^{r_{2}t}-\frac{1}{2\nu|\xi|^{2}}
\mathrm{e}^{r_{3}t}-\frac{1}{2\nu|\xi|^{2}}\mathrm{e}^{r_{4}t}\Big)\hat{\rho}_{0}^{h}\\ &+\Big(\frac{\alpha_{1}}{\nu^{2}|\xi|^{4}}(\alpha_{2}+\frac{\alpha_1}{\nu})\mathrm{e}^{r_{1}t}-\frac{1}{\nu|\xi|^{2}}\mathrm{e}^{r_{2}t}
+i\frac{1}{2\nu|\xi|^{3}}\mathrm{e}^{r_{3}t}-i\frac{1}{2\nu|\xi|^{3}}\mathrm{e}^{r_{4}t}\Big)
i\frac{\xi \hat{m}_{0}^{h}}{|\xi |} \\
&+\Big(\frac{\alpha_{1}}{\nu|\xi|}\mathrm{e}^{r_{1}t}-\frac{\alpha_{1}}{\nu|\xi|}\mathrm{e}^{r_{2}t}-\frac{\alpha_{1}\alpha_{2}}{2\nu|\xi|}
\mathrm{e}^{r_{3}t}-\frac{\alpha_{1}\alpha_{2}}{2\nu|\xi|}\mathrm{e}^{r_{4}t}\Big)\hat{n}_{0}^{h}\\ &+\Big(-\frac{\alpha_{2}(1+\nu+\alpha_{1})}{\nu|\xi|^{2}}\mathrm{e}^{r_{1}t}+\mathrm{e}^{r_{2}t}-\frac{1}{8|\xi|^{2}}\mathrm{e}^{r_{3}t}
-\frac{1}{8|\xi|^{2}}\mathrm{e}^{r_{4}t}\Big)i\frac{\xi \hat{\omega}_{0}^{h}}{|\xi |}+\cdots.\ \ \ \ \ \ \ \ \ \
\end{aligned}
\end{equation}
\end{lemma}

\subsection{Spectral analysis for the incompressible part}

We can express the IVP  (\ref{2.8})  for  $\mathcal{V}=(\Phi, \Psi)^{T}$  as
\begin{equation}\label{2.26}
\left\{\begin{array}{l}
\mathcal{V}_{t}=\mathcal{A}_{2} \mathcal{V}, \\
\left.\mathcal{V}\right|_{t=0}=\mathcal{V}_{0},
\end{array}\right.
\end{equation}
where the operator  $\mathcal{A}_{1}$  is given by
\begin{equation}
\mathcal{A}_{2}=\left(\begin{array}{cc}
-1 & \alpha_{2} \\
1 & -\alpha_{2}-\bar{\mu} \Lambda^{2}
\end{array}\right).
\end{equation}
Applying the Fourier transform to the system  (\ref{2.26}), one has
\begin{equation}
\left\{\begin{array}{l}
\hat{\mathcal{V}}_{t}=\mathcal{A}_{2}(\xi) \hat{\mathcal{V}}, \\
\left.\hat{\mathcal{V}}\right|_{t=0}=\hat{\mathcal{V}}_{0},
\end{array}\right.
\end{equation}
where $\mathcal{A}_{2}(\xi)$ is defined by
\begin{equation}
\mathcal{A}_{2}(\xi)=\left(\begin{array}{cc}
-1 & \alpha_{2} \\
1 & -\alpha_{2}-\bar{\mu}|\xi|^{2}
\end{array}\right),
\end{equation}
and its eigenvalues satisfy
\begin{equation}
\operatorname{det}\left(\kappa I-\mathcal{A}_{2}(\xi)\right)=\kappa^{2}+\left(\alpha_{2}+1+\bar{\mu}|\xi|^{2}\right) \kappa+\bar{\mu}|\xi|^{2}=0.
\end{equation}
Then, one can decompose the semigroup $e^{t\mathcal{A}_{2}(\xi)}$ as follows:
\begin{equation}\label{2.31}
e^{t\mathcal{A}_{2}(\xi)}=e^{\kappa_1t}Q_1+e^{\kappa_2t}Q_2,
\end{equation}
with
\begin{equation}\label{2.32}
Q_1(\xi)=\frac{\mathcal{A}_{2}(\xi)-\kappa_2I}{\kappa_1-\kappa_2},\ \ Q_2(\xi)=\frac{\mathcal{A}_{2}(\xi)-\kappa_1I}{\kappa_2-\kappa_1}.
\end{equation}

By a direct computation, one can have the following two lemmas on the asymptotic expansion in low frequency and high frequency parts respectively.
\begin{lemma}\label{l 2.6}
 There exists a positive constant  $\eta_{1} \ll 1$  such that, for  $|\xi| \ll \eta_{1}$, the spectral has the following Taylor series expansion:
\begin{equation}\label{2.37}
\left\{\begin{array}{ll}
\kappa_{1}=-\alpha_{2}-1-\frac{\alpha_{2} \bar{\mu}}{\alpha_{2}+1}|\xi|^{2}+O\left(|\xi|^{4}\right), \\[2mm]
\kappa_{2}=-\frac{\bar{\mu}}{\alpha_{2}+1}|\xi|^{2}+O\left(|\xi|^{4}\right),
\end{array}\right.
\end{equation}
\begin{equation*}
Q_{1}^l(\xi)=\frac{\mathcal{A}_{2}(\xi)-\kappa_{2} I}{\kappa_{1}-\kappa_{2}}=-\frac{1}{\left(\alpha_{2}+1\right)}\left(\begin{array}{cc}
-1 & \alpha_{2} \\
 1& -\alpha_{2}
\end{array}\right)+\mathcal{O}(|\xi|),
\end{equation*}
\begin{equation*}
Q_{2}^l(\xi)=\frac{\mathcal{A}_{2}(\xi)-\kappa_{1} I}{\kappa_{2}-\kappa_{1}}=\frac{1}{\left(\alpha_{2}+1\right)}\left(\begin{array}{ll}
\alpha_{2} & \alpha_{2} \\
1 & 1
\end{array}\right)+\mathcal{O}(|\xi|),
\end{equation*}
\begin{equation*}
\begin{aligned}
\big(\hat{\Phi}^{jk}\big)^l=&\frac{\mathrm{e}^{\kappa_{1}t}}{\alpha_{2}+1}\big( \hat{\Phi}_{0}^{jk}-\alpha_{2} \hat{\Psi}_{0}^{jk}\big)+\frac{\alpha_{2}\mathrm{e}^{\kappa_{2}t}}{\alpha_{2}+1}\big( \hat{\Phi}_{0}^{jk}+\hat{\Psi}_{0}^{jk}\big)+\cdots\\[2mm]
=&\frac{1}{\alpha_{2}+1}( \mathrm{e}^{\kappa_{1}t}+\alpha_{2}\mathrm{e}^{\kappa_{2}t })i\frac{\xi^{k}u_{0}^{j}-\xi^{j}u_{0}^{k}}{|\xi|}+\frac{\alpha_{2}}{\alpha_{2}+1}( -\mathrm{e}^{\kappa_{1}t}+\mathrm{e}^{\kappa_{2}t})i\frac{\xi^{k}v_{0}^{j}-\xi^{j}v_{0}^{k}}{|\xi|}+\cdots,
\end{aligned}
\end{equation*}
\begin{equation*}
\begin{aligned}
\big(\hat{\Psi}^{jk}\big)^l=&\frac{\mathrm{e}^{\kappa_{1}t}}{\alpha_{2}+1}( -\hat{\Phi}_{0}^{jk}+ \alpha_{2}\hat{\Psi}_{0}^{jk})+\frac{\mathrm{e}^{\kappa_{2}t}}{\alpha_{2}+1}( \hat{\Phi}_{0}^{jk}+\hat{\Psi}_{0}^{jk})+\cdots\\[2mm]
=&\frac{\alpha_{2}}{\alpha_{2}+1}( -\mathrm{e}^{\kappa_{1}t}+\mathrm{e}^{\kappa_{2}t })i\frac{\xi^{k}u_{0}^{j}-\xi^{j}u_{0}^{k}}{|\xi|}+\frac{1}{\alpha_{2}+1}( \alpha_{2}\mathrm{e}^{\kappa_{1}t}+\mathrm{e}^{\kappa_{2}t})i\frac{\xi^{k}v_{0}^{j}-\xi^{j}v_{0}^{k}}{|\xi|}+\cdots.
\end{aligned}
\end{equation*}
\end{lemma}

\begin{lemma}\label{l 2.7}
 There exists a positive constant  $\eta_{2} \gg 1$  such that, for  $|\xi| \gg \eta_{2}$, the spectral has the following Taylor series expansion:
\begin{equation}\label{2.38}
\left\{\begin{array}{l}
\kappa_{1}=-1+O\left(|\xi|^{-2}\right),\\[2mm]
\kappa_{2}=-\bar{\mu}|\xi|^{2}-\alpha_{2}+O\left(|\xi|^{-2}\right).
\end{array}\right.
\end{equation}
\begin{equation*}
Q_{1}^h(\xi)=\frac{\mathcal{A}_{2}(\xi)-\kappa_{2} I}{\kappa_{1}-\kappa_{2}}=\frac{1}{-1+\alpha_{2}+\bar{\mu}|\xi|^{2}}\left(\begin{array}{cc}
\bar{\mu}|\xi|^{2} & \alpha_{2} \\
1 & 0
\end{array}\right)+\mathcal{O}(|\xi|^{-2}),
\end{equation*}
\begin{equation*}
Q_{2}^h(\xi)=\frac{\mathcal{A}_{2}(\xi)-\kappa_{1} I}{\kappa_{2}-\kappa_{1}}=\frac{-1}{-1+\alpha_{2}+\bar{\mu}|\xi|^{2}}\left(\begin{array}{ll}
0 & \ \ \ \ \alpha_{2} \\
1 & -\bar{\mu}|\xi|^{2}
\end{array}\right)+\mathcal{O}(|\xi|^{-2}),
\end{equation*}
\begin{equation}
\begin{aligned}
\big(\hat{\Phi}^{jk}\big)^h=&\frac{\mathrm{e}^{\kappa_{1}t}}{-1+\alpha_{2}+\bar{\mu}|\xi|^{2}}(\bar{\mu}|\xi|^{2}\hat{\Phi}^{jk}_{0}+ \alpha_{2}\hat{\Psi}^{jk}_{0})-\frac{\mathrm{e}^{\kappa_{2}t}}{-1+\alpha_{2}+\bar{\mu}|\xi|^{2}}( \alpha_{2}\hat{\Psi}^{jk}_{0})+\cdots\\[2mm]
=&\frac{\bar{\mu}|\xi|^{2}\mathrm{e}^{\kappa_{1}t}}{-1+\alpha_{2}+\bar{\mu}|\xi|^{2}}i\frac{\xi^{k}u_{0}^{j}-\xi^{i}u_{0}^{j}}{|\xi|}
+\frac{\alpha_{2}\mathrm{e}^{\kappa_{1}t}-\alpha_{2}\mathrm{e}^{\kappa_{2}t}}{-1+\alpha_{2}+\bar{\mu}|\xi|^{2}}
i\frac{\xi^{k}v_{0}^{j}-\xi^{j}v_{0}^{k}}{|\xi|}+\cdots,
\end{aligned}
\end{equation}
\begin{equation}
\begin{aligned}
\big(\hat{\Psi}^{jk}\big)^h=&\frac{\mathrm{e}^{\kappa_{1}t}}{-1+\alpha_{2}+\bar{\mu}|\xi|^{2}}\hat{\Phi}_{0}^{jk}
-\frac{\mathrm{e}^{\kappa_{2}t}}{-1+\alpha_{2}+\bar{\mu}|\xi|^{2}}( \hat{\Phi}_{0}^{jk}-\bar{\mu}|\xi|^{2}\hat{\Psi}_{0}^{jk})+\cdots\\[2mm]
=&\frac{\mathrm{e}^{\kappa_{1}t}-\mathrm{e}^{\kappa_{2}t }}{-1+\alpha_{2}+\bar{\mu}|\xi|^{2}}i\frac{\xi^{k}u_{0}^{j}-\xi^{j}u_{0}^{k}}{|\xi|}
+\frac{\bar{\mu}|\xi|^{2}\mathrm{e}^{\kappa_{2}t}}{-1+\alpha_{2}+\bar{\mu}|\xi|^{2}} i\frac{\xi^{k}v_{0}^{j}-\xi^{j}v_{0}^{k}}{|\xi|}+\cdots.
\end{aligned}
\end{equation}
\end{lemma}

With the spectral analysis in hand, we can derive the pointwise space-time estimates for each component in the Green's function.

\section{Pointwise space-time behavior of Green's function}

\subsection{Representation of the solution in Fourier space}

The representations of two densities $\hat{\rho}$ and $\hat{n}$ in the low frequency and high frequency parts have been given in (\ref{2.21})-(\ref{2.22}) and (\ref{2.27})-(\ref{2.28}). Foe two momenta $\hat{m}$ and $\hat{\omega}$, due to the Hodge decomposition, we can get them as follows from (\ref{2.23})-(\ref{2.24}) and (\ref{2.29})-(\ref{2.30}):
\begin{equation}\label{3.1}
\begin{aligned}
\hat{m}^l =&-\widehat{\wedge^{-1}\nabla \varphi}-\widehat{\wedge ^{-1}{\rm div} \Phi}\\
=&-\frac{1}{2\left(\alpha_{1}+\alpha_{2}\right)}\left(-i \alpha_{2} \sqrt{\frac{\alpha_{1}+\alpha_{2}}{\alpha_{2}+1}}\mathrm{e}^{r_{3}t}+i \alpha_{2} \sqrt{\frac{\alpha_{1} +\alpha_{2}}{\alpha_{2}+1}}\mathrm{e}^{r_{4} t} \right)\frac{i\xi\hat{\rho}_{0}^{l}}{|\xi|}\\
&+\frac{1}{2\left(\alpha_{2}+1\right)}\left(2\mathrm{e}^{r_{1} t}+\alpha_{2} \mathrm{e}^{r_{3}t}+\alpha_{2} \mathrm{e}^{r_{4}t}\right)\frac{\xi\xi^T\hat{m}_{0}^{l}}{|\xi|^{2}}\\
&-\frac{1}{2\left(\alpha_{1}+\alpha_{2}\right)}\left(-i \alpha_{1}\alpha_{2} \sqrt{\frac{\alpha_{1}+\alpha_{2}}{\alpha_{2}+1}}\mathrm{e}^{r_{3}t}+i \alpha_{1} \alpha_{2}\sqrt{\frac{\alpha_{1} +\alpha_{2}}{\alpha_{2}+1}}\mathrm{e}^{r_{4} t} \right)\frac{i\xi\hat{n}_{0}^{l}}{|\xi|}\\
&+\frac{1}{2\left(\alpha_{2}+1\right)}\left(-2\alpha_{2}\mathrm{e}^{r_{1} t}+\alpha_{2}\mathrm{e}^{r_{3}t}+\alpha_{2} \mathrm{e}^{r_{4}t}\right)\frac{\xi\xi^T \hat{\omega}_{0}^{l}}{|\xi|^{2}}\\
&+\frac{1}{\alpha_{2}+1}\left( \mathrm{e}^{\kappa_{1}t}+\alpha_{2}\mathrm{e}^{\kappa_{2}t }\right)\Big(I-\frac{\xi\xi^T}{|\xi|^{2}}\Big)\hat{m}_{0}^{l}+\frac{\alpha_{2}}{\alpha_{2}+1}\left( -\mathrm{e}^{\kappa_{1}t}+\mathrm{e}^{\kappa_{2}t}\right)\Big(I-\frac{\xi\xi^T}{|\xi|^{2}}\Big)\hat{\omega}_{0}^{l}+\cdots,
\end{aligned}
\end{equation}
\begin{equation}\label{3.2}
\begin{aligned}
\hat{\omega}^l =&-\widehat{\wedge^{-1}\nabla \psi}-\widehat{\wedge ^{-1}{\rm div} \Psi}\\
=&-\frac{1}{2\left(\alpha_{1}+\alpha_{2}\right)}\left(-i  \sqrt{\frac{\alpha_{1}+\alpha_{2}}{\alpha_{2}+1}}\mathrm{e}^{r_{3}t}+i  \sqrt{\frac{\alpha_{1} +\alpha_{2}}{\alpha_{2}+1}}\mathrm{e}^{r_{4} t}\right)\frac{i\xi\hat{\rho}_{0}^{l}}{|\xi|}\\
&+\frac{1}{2\left(\alpha_{2}+1\right)}\left(-2\mathrm{e}^{r_{1} t}+ \mathrm{e}^{r_{3}t}+\mathrm{e}^{r_{4}t}\right)\frac{\xi\xi^T\hat{m}_{0}^{l}}{|\xi|^{2}}\\
&-\frac{1}{2\left(\alpha_{1}+\alpha_{2}\right)}\left(-i \alpha_{1} \sqrt{\frac{\alpha_{1}+\alpha_{2}}{\alpha_{2}+1}}\mathrm{e}^{r_{3}t}+i \alpha_{1} \sqrt{\frac{\alpha_{1} +\alpha_{2}}{\alpha_{2}+1}}\mathrm{e}^{r_{4} t}\right)\frac{i\xi\hat n_{0}^{l}}{|\xi|}\\
&+\frac{1}{2\left(\alpha_{2}+1\right)}\left(2\alpha_{2}\mathrm{e}^{r_{1} t}+\mathrm{e}^{r_{3}t}+ \mathrm{e}^{r_{4}t}\right)\frac{\xi\xi^T \hat{\omega}_{0}^{l}}{|\xi|^{2}}\\
&+\frac{\alpha_{2}}{\alpha_{2}+1}\left( -\mathrm{e}^{\kappa_{1}t}+\mathrm{e}^{\kappa_{2}t }\right)\Big(I-\frac{\xi\xi^T}{|\xi|^{2}}\Big)\hat{m}_{0}^{l}+\frac{1}{\alpha_{2}+1}\left( \alpha_{2}\mathrm{e}^{\kappa_{1}t}+\mathrm{e}^{\kappa_{2}t}\right)\Big(I-\frac{\xi\xi^T}{|\xi|^{2}}\Big)\hat{\omega}_{0}^{l}+\cdots,
\end{aligned}
\end{equation}
\begin{equation}\label{3.3}
\begin{aligned}
\hat{m}^h =&-\widehat{\wedge^{-1}\nabla \varphi}-\widehat{\wedge ^{-1}{\rm div} \Phi}\\
=&-\left(\frac{\alpha_1^{2}}{\nu^{3}|\xi|^{3}}\mathrm{e}^{r_{1}t}-\frac{\alpha_1}{\nu^{4}|\xi|^{5}}\left(\frac{\alpha_{1}}{\nu}-\alpha_{2}\right)\mathrm{e}^{r_{2}t} -i\frac{1}{2}\mathrm{e}^{r_{3} t} +i\frac{1}{2}\mathrm{e}^{r_{4} t}\right)\frac{i\xi\hat{\rho}_{0}^{h}}{|\xi|}\\
&+\left(\frac{1}{\nu^{2}|\xi|^{2}}\big(\frac{\alpha_{1}}{\nu}-\alpha_{2}\big)\mathrm{e}^{r_{1}t}+\frac{\alpha_{2}}{\nu^{2}|\xi|^{4}}\mathrm{e}^{r_{2}t} +\frac{1}{2}\mathrm{e}^{r_{3} t} +\frac{1}{2}\mathrm{e}^{r_{4} t}\right)\frac{\xi\xi^T\hat{m}_{0}^{h}}{|\xi|^{2}}\\
&-\left(-\frac{\alpha_1^{2}\alpha_{2}}{\nu^{2}|\xi|^{3}}\mathrm{e}^{r_{1}t}+\frac{\alpha_{1}\alpha_{2}}{\nu^{2}|\xi|^{3}}\mathrm{e}^{r_{2}t} -i\frac{\alpha_{1}\alpha_{2}}{2\nu|\xi|^{2}}\mathrm{e}^{r_{3}t}+i\frac{\alpha_{1}\alpha_{2}}{2\nu|\xi|^{2}}\mathrm{e}^{r_{4}t} \right)\frac{i\xi\hat{n}_{0}^{h}}{|\xi|}\\
&+\left(\Big(\frac{\alpha_{1}\alpha_{2}}{\nu}+(\frac{\alpha_{1}}{\nu})^{2}-1-\alpha_{2}\Big)\frac{\alpha_{2}\mathrm{e}^{r_{1}t}}{\nu|\xi|^{4}}
-\frac{\alpha_{2}\mathrm{e}^{r_{2}t}}{\nu|\xi|^{2}} +\frac{\alpha_{2}\mathrm{e}^{r_{3}t}}{2\nu|\xi|^{2}}+\frac{\mathrm{e}^{r_{4}t}}{2\nu|\xi|^{2}} \right)\frac{\xi\xi^T \hat{\omega}_{0}^{h}}{|\xi|^{2}}\\
&+\frac{\bar{\mu}|\xi|^{2}\mathrm{e}^{\kappa_{1}t}}{-1+\alpha_{2}+\bar{\mu}|\xi|^{2}}\Big(I-\frac{\xi\xi^T}{|\xi|^{2}}\Big)\hat{m}_{0}^{h}
+\frac{\alpha_{2}\mathrm{e}^{\kappa_{1}t}-\alpha_{2}\mathrm{e}^{\kappa_{2}t}}{-1+\alpha_{2}
+\bar{\mu}|\xi|^{2}}\Big(I-\frac{\xi\xi^T}{|\xi|^{2}}\Big)\hat{\omega}_{0}^{h}+\cdots,\ \ \ \ \ \ \ \ \ \
\end{aligned}
\end{equation}
\begin{equation}\label{3.4}
\begin{aligned}
\hat{\omega}^h =&-\widehat{\wedge^{-1}\nabla \psi}-\widehat{\wedge ^{-1}{\rm div} \Psi}\\
=&-\left(-\frac{\alpha_{1}}{\nu^{2}|\xi|^{3}}\mathrm{e}^{r_{1}t}
+\frac{1}{\nu^{2}|\xi|^{3}}\mathrm{e}^{r_{2}t}-\frac{1}{2\nu|\xi|^{2}}\mathrm{e}^{r_{3}t}
-\frac{1}{2\nu|\xi|^{2}}\mathrm{e}^{r_{4}t}\right)\frac{i\xi\hat{\rho}_{0}^{h}}{|\xi|}\\
&+\left(\frac{\alpha_{1}}{\nu^{2}|\xi|^{4}}(\alpha_{2}+\frac{\alpha_1}{\nu})\mathrm{e}^{r_{1}t}-\frac{1}{\nu|\xi|^{2}}\mathrm{e}^{r_{2}t}
+i\frac{1}{2\nu|\xi|^{3}}\mathrm{e}^{r_{3}t}-i\frac{1}{2\nu|\xi|^{3}}\mathrm{e}^{r_{4}t}\right)\frac{\xi\xi^T\hat{m}_{0}^{h}}{|\xi|^{2}}\\
&-\left(\frac{\alpha_{1}}{\nu|\xi|}\mathrm{e}^{r_{1}t}-\frac{\alpha_{1}}{\nu|\xi|}\mathrm{e}^{r_{2}t}
-\frac{\alpha_{1}\alpha_{2}}{2\nu|\xi|}\mathrm{e}^{r_{3}t}-\frac{\alpha_{1}\alpha_{2}}{2\nu|\xi|}\mathrm{e}^{r_{4}t}\right)\frac{i\xi\hat{n}_{0}^{h}}{|\xi|}\\ &+\left(-\frac{\alpha_{2}(1+\nu+\alpha_{1})}{\nu|\xi|^{2}}\mathrm{e}^{r_{1}t}
+\mathrm{e}^{r_{2}t}-\frac{1}{8|\xi|^{2}}\mathrm{e}^{r_{3}t}-\frac{1}{8|\xi|^{2}}\mathrm{e}^{r_{4}t}\right)\frac{\xi\xi^T \hat{\omega}_{0}^{h}}{|\xi|^{2}}\\
&+\frac{ \mathrm{e}^{\kappa_{1}t}-\mathrm{e}^{\kappa_{2}t }}{-1+\alpha_{2}+\bar{\mu}|\xi|^{2}}\Big(I-\frac{\xi\xi^T}{|\xi|^{2}}\Big)\hat{m}_{0}^{h}
+\frac{\bar{\mu}|\xi|^{2}\mathrm{e}^{\kappa_{2}t}}{-1+\alpha_{2}
+\bar{\mu}|\xi|^{2}}\Big(I-\frac{\xi\xi^T}{|\xi|^{2}}\Big)\hat{\omega}_{0}^{h}+\cdots.\ \ \ \ \ \ \ \ \ \
\end{aligned}
\end{equation}

\subsection{Pointwise space-time estimate of Green's matrix}

In this section, we will give the pointwise estimate of each component in Green's matrix by using the spectrum analysis in Section 2. The difficulties mainly include using suitable combination to avoid the singularity in the low frequency part arising from the Hodge decomposition, and giving the description of the singular part in the high frequency arising from the definition of the Green's function. For convenience in the proof of the pointwise space-time estimates for the nonlinear problem, we shall first give the following estimates for the Green's matrix in the Fourier space.

\bigbreak

\textbf{Pointwise description of low frequency part.}
\vspace{5mm}\\
From (\ref{2.21})-(\ref{2.22}), we have
\begin{equation}\label{3.5}
\begin{array}{rl}
&\hat{G}_{11}^l\sim \frac{ 2\alpha_{1}\mathrm{e}^{r_{2} t}+\alpha_{2}\mathrm{e}^{r_{3} t}+\alpha_{2}\mathrm{e}^{r_{4} t}}{2(\alpha_{1}+\alpha_{2})},\ \ \
\hat{G}_{12}^l\sim \frac{-\alpha_{2}}{2(\alpha_{1}+\alpha_{2})}\sqrt{\frac{\alpha_{1} +\alpha_{2}}{\alpha_{2}+1}}( \mathrm{e}^{r_{3} t}-\mathrm{e}^{r_{4} t} )\frac{\xi }{|\xi |},\\[2mm]
&\hat{G}_{13}^l\sim \frac{\alpha_{1}\alpha_{2}}{2(\alpha_{1}+\alpha_{2})}\big(-2\mathrm{e}^{r_{2} t}+\mathrm{e}^{r_{3} t}+\mathrm{e}^{r_{4} t} \big),\ \
\hat{G}_{14}^l\sim \frac{-\alpha_{2}}{2(\alpha_{1}+\alpha_{2})}\sqrt{\frac{\alpha_{1} +\alpha_{2}}{\alpha_{2}+1}}(\mathrm{e}^{r_{3} t}-\mathrm{e}^{r_{4} t})\frac{\xi }{|\xi |},
\end{array}
\end{equation}
and
\begin{equation}\label{3.6}
\begin{array}{rl}
&\hat{G}_{31}^l\sim \frac{ -2\mathrm{e}^{r_{2} t}+\mathrm{e}^{r_{3} t}+\mathrm{e}^{r_{4} t}}{2(\alpha_{1}+\alpha_{2})},\ \ \
\hat{G}_{32}^l\sim \frac{-1}{2(\alpha_{1}+\alpha_{2})}\sqrt{\frac{\alpha_{1} +\alpha_{2}}{\alpha_{2}+1}}( \mathrm{e}^{r_{3} t}-\mathrm{e}^{r_{4} t} )\frac{\xi }{|\xi |},\\[2mm]
&\hat{G}_{33}^l\sim \frac{1}{2(\alpha_{1}+\alpha_{2})}\big(2\alpha_2\mathrm{e}^{r_{2} t}+\alpha_1\mathrm{e}^{r_{3} t}+\alpha_1\mathrm{e}^{r_{4} t} \big),\ \ \
\hat{G}_{34}^l\sim \frac{-1}{2(\alpha_{1}+\alpha_{2})}\sqrt{\frac{\alpha_{1} +\alpha_{2}}{\alpha_{2}+1}}(\mathrm{e}^{r_{3} t}-\mathrm{e}^{r_{4} t})\frac{\xi }{|\xi |}.
\end{array}
\end{equation}
In the same way, we can give the estimate for two momenta from (\ref{3.1})-(\ref{3.2}):
\begin{equation}\label{3.7}
\begin{array}{rl}
&\hat{G}_{21}^l\sim \frac{-\alpha_{2}}{2\left(\alpha_{1}+\alpha_{2}\right)}\sqrt{\frac{\alpha_{1}+\alpha_{2}}{\alpha_{2}+1}}\left( \mathrm{e}^{r_{3}t}-\mathrm{e}^{r_{4} t} \right)\frac{\xi}{|\xi|},\\
&\hat{G}_{22}^l\sim \frac{1}{2\left(\alpha_{2}+1\right)}\left(2\mathrm{e}^{r_{1} t}+\alpha_{2} \mathrm{e}^{r_{3}t}+\alpha_{2} \mathrm{e}^{r_{4}t}\right)\frac{\xi\xi^T}{|\xi|^{2}}+\frac{\mathrm{e}^{\kappa_{1}t}+\alpha_{2}\mathrm{e}^{\kappa_{2}t }}{\alpha_{2}+1}\Big(I-\frac{\xi\xi^T}{|\xi|^{2}}\Big),\\[2mm]
&\ \ \ \ \ \ \ \sim \frac{\alpha_2}{\alpha_{2}+1}\big(\frac{\mathrm{e}^{r_{3}t}+ \mathrm{e}^{r_{4}t}}{2}-e^{\kappa_2t}\big)\frac{\xi\xi^T}{|\xi|^{2}}+\frac{1}{\alpha_2+1}
(\mathrm{e}^{r_1t}-\mathrm{e}^{\kappa_1t})\frac{\xi\xi^T}{|\xi|^{2}}+\frac{\alpha_2}{\alpha_2+1}\mathrm{e}^{\kappa_2t}
+\frac{1}{\alpha_2+1}\mathrm{e}^{\kappa_1t},\\[2mm]
&\hat{G}_{23}^l\sim \frac{-\alpha_{1}\alpha_{2} }{2\left(\alpha_{1}+\alpha_{2}\right)}\sqrt{\frac{\alpha_{1}+\alpha_{2}}{\alpha_{2}+1}}\left( \mathrm{e}^{r_{3}t}-\mathrm{e}^{r_{4} t} \right)\frac{\xi}{|\xi|},\\[2mm]
&\hat{G}_{24}^l\sim \frac{\alpha_{2}}{2\left(\alpha_{2}+1\right)}\left(-2\mathrm{e}^{r_{1} t}+\mathrm{e}^{r_{3}t}+ \mathrm{e}^{r_{4}t}\right)\frac{\xi\xi^T}{|\xi|^{2}}+\frac{\alpha_{2}}{\alpha_{2}+1}\left( -\mathrm{e}^{\kappa_{1}t}+\mathrm{e}^{\kappa_{2}t}\right)\Big(I-\frac{\xi\xi^T}{|\xi|^{2}}\Big),\\[2mm]
&\ \ \ \ \ \ \ \sim \frac{\alpha_2}{\alpha_{2}+1}\big(\frac{\mathrm{e}^{r_{3}t}+ \mathrm{e}^{r_{4}t}}{2}-\mathrm{e}^{\kappa_2t}\big)\frac{\xi\xi^T}{|\xi|^{2}}
+\frac{-\alpha_2}{\alpha_2+1}(\mathrm{e}^{r_1t}-\mathrm{e}^{\kappa_1t})\frac{\xi\xi^T}{|\xi|^{2}}
-\frac{\alpha_2}{\alpha_2+1}(\mathrm{e}^{\kappa_1t}-\mathrm{e}^{\kappa_2t}),
\end{array}
\end{equation}
and
\begin{equation}\label{3.8}
\begin{array}{rl}
&\hat{G}_{41}^l\sim \frac{-1}{2\left(\alpha_{1}+\alpha_{2}\right)}\sqrt{\frac{\alpha_{1}+\alpha_{2}}{\alpha_{2}+1}}\left( \mathrm{e}^{r_{3}t}-\mathrm{e}^{r_{4} t} \right)\frac{\xi}{|\xi|},\\
&\hat{G}_{42}^l\sim \frac{1}{2\left(\alpha_{2}+1\right)}\left(-2\mathrm{e}^{r_{1} t}+\mathrm{e}^{r_{3}t}+ \mathrm{e}^{r_{4}t}\right)\frac{\xi\xi^T}{|\xi|^{2}}+\frac{-\alpha_{2}\mathrm{e}^{\kappa_{1}t}+\alpha_{2}\mathrm{e}^{\kappa_{2}t }}{\alpha_{2}+1}\Big(I-\frac{\xi\xi^T}{|\xi|^{2}}\Big),\\[2mm]
&\ \ \ \ \ \ \ \sim \frac{1}{\alpha_{2}+1}\big(\frac{\mathrm{e}^{r_{3}t}+ \mathrm{e}^{r_{4}t}}{2}-e^{\kappa_2t}\big)\frac{\xi\xi^T}{|\xi|^{2}}+\frac{1}{\alpha_2+1}(\mathrm{e}^{\kappa_{2}t}
-\mathrm{e}^{r_{1}t})+\frac{-\alpha_{2}(\mathrm{e}^{\kappa_{1}t}
-\mathrm{e}^{\kappa_{2}t})}{\alpha_{2}+1}\Big(I-\frac{\xi\xi^T}{|\xi|^{2}}\Big),\\[2mm]
&\hat{G}_{43}^l\sim \frac{-\alpha_{1}}{2\left(\alpha_{1}+\alpha_{2}\right)}\sqrt{\frac{\alpha_{1}+\alpha_{2}}{\alpha_{2}+1}}\left( \mathrm{e}^{r_{3}t}-\mathrm{e}^{r_{4} t} \right)\frac{\xi}{|\xi|},\\[2mm]
&\hat{G}_{44}^l\sim \frac{1}{2\left(\alpha_{2}+1\right)}\left(2\alpha_{2}\mathrm{e}^{r_{1} t}+\mathrm{e}^{r_{3}t}+ \mathrm{e}^{r_{4}t}\right)\frac{\xi\xi^T}{|\xi|^{2}}+\frac{1}{\alpha_{2}+1}\left( \alpha_{2}\mathrm{e}^{\kappa_{1}t}+\mathrm{e}^{\kappa_{2}t}\right)\Big(I-\frac{\xi\xi^T}{|\xi|^{2}}\Big),\\[2mm]
&\ \ \ \ \ \ \ \sim \frac{1}{\alpha_{2}+1}\big(\frac{\mathrm{e}^{r_{3}t}+ \mathrm{e}^{r_{4}t}}{2}-e^{\kappa_2t}\big)\frac{\xi\xi^T}{|\xi|^{2}}
+\frac{\alpha_2}{\alpha_2+1}(e^{r_1t}-e^{\kappa_1t})\frac{\xi\xi^T}{|\xi|^{2}}+\frac{\alpha_2}{\alpha_2+1}\mathrm{e}^{\kappa_1t}
+\frac{1}{\alpha_2+1}\mathrm{e}^{\kappa_2t},
\end{array}
\end{equation}
where $``\sim"$ means that only the leading terms are stated at the right hand side, since the rest terms don't affect the results. After these recombinations, the singularity at $\xi=0$ in the above terms  can be treated, which is crucial for us to derive the pointwise space-time estimates for the Green's function. We only take several terms for examples. The first one is $( \mathrm{e}^{r_{3} t}-\mathrm{e}^{r_{4} t} )\frac{\xi }{|\xi |}$ in $\hat{G}_{12}^l$, which together with (\ref{2.16}) and Euler formula implies that $( \mathrm{e}^{r_{3} t}-\mathrm{e}^{r_{4} t} )\frac{\xi }{|\xi |}\sim \sin (c|\xi|t)\frac{\xi }{|\xi |}=\frac{\sin(c|\xi|t}{c|\xi|}\cdot c\xi:=c\xi\cdot\hat{\mathbf{w}}$. Here and below, we denote $c=\sqrt{\frac{\alpha_{1}+\alpha_{2}}{\alpha_{2}+1}}$, which can be regarded as the basic sound speed for this two-phase model. The second one is $\big(\frac{\mathrm{e}^{r_{3}t}+ \mathrm{e}^{r_{4}t}}{2}-e^{\kappa_2t}\big)\frac{\xi\xi^T}{|\xi|^{2}}$ in $\hat{G}_{22}^l$, which together with the asymptotic expansions  in (\ref{2.16}) and (\ref{2.37}) for $r_3,r_4,\kappa_2$ and Euler formula implies that it can be rewritten as $\cos(c|\xi|t)\mathrm{e}^{-|\xi|^2t}\frac{\xi\xi^T}{|\xi|^{2}}$. Then it will be solved by splitting it into Riesz wave-$I$ with the symbol $(\cos(c|\xi| t)-1)\frac{\xi\xi^T}{|\xi|^2}e^{-\theta_1|\xi|^2t}$ and Riesz wave-$I\!I$ with the symbol $\big(e^{-\theta_1|\xi|^2t}-e^{-\theta_2|\xi|^2t}\big)\ \frac{\xi\xi^T}{|\xi|^2}$. The readers can see the details in Lemma 4.7 and Lemma 4.8 of Du-Wu\cite{Du}, which yields that its pointwise space-time description contains both the Huygens' wave and the Riesz wave (diffusion wave).

Thus, we can get the pointwise estimate for the Green's function in the low frequency:
\begin{lemma}\label{l 3.1} For any $|\alpha|\geq0$, there exists a constant $C>0$ such that
\begin{equation*}
\begin{array}{rl}
&|D_x^\alpha(G_{11}^l,G_{13}^l,G_{24}^l, G_{31}^l,G_{33}^l,G_{44}^l)|\\[2mm]
&\ \ \ \ \ \ \ \ \leq C(1+t)^{-\frac{3+|\alpha|}{2}}\Big(1+\frac{|x|^2}{1+t}\Big)^{-N}+C(1+t)^{-\frac{4+|\alpha|}{2}}\Big(1+\frac{(|x|-ct)^2}{1+t}\Big)^{-N},\\[2mm]
&|D_x^\alpha (G_{22}^l,G_{42}^l)|\leq C(1+t)^{-\frac{3+|\alpha|}{2}}\Big(1+\frac{|x|^2}{1+t}\Big)^{-\frac{3+|\alpha|}{2}}+C(1+t)^{-\frac{4+|\alpha|}{2}}\Big(1+\frac{(|x|-ct)^2}{1+t}\Big)^{-N},\ \ \ \\[2mm]
&|D_x^\alpha(G_{12}^l,G_{14}^l,G_{21}^l,G_{23}^l, G_{32}^l,G_{34}^l,G_{41}^l,G_{43}^l)|\leq C(1+t)^{-\frac{4+|\alpha|}{2}}\Big(1+\frac{(|x|-ct)^2}{1+t}\Big)^{-N},
\end{array}
\end{equation*}
where the constant $N>0$ can be arbitrarily large.
\end{lemma}

However, the estimates in Lemma \ref{l 3.1} for the low frequency part  don't suffices to yield the desired pointwise space-time description, since the non-conservation of the system (\ref{1.1}). In fact, the conservation is crucial when dealing with nonlinear convolution containing the Huygens' wave, see the nonlinear  estimate $K_3$ in Lemma \ref{A.5}. Thus, we have to seek some additional decay rates from the corresponding components in the Green's function.

 We find that the non-conservation is just from two damping terms in $F_1$ and $F_2$ of the momentum equations in (\ref{2.1})-(\ref{2.2}), and fortunately they have the opposite sign. This gives us a little hope to get extra decay from the low frequency part. To this end, we recall the second and the forth columns of the Green's function containing $\mathrm{e}^{r_{3}t}$ and $\mathrm{e}^{r_{4}t}$ in (\ref{3.5})-(\ref{3.8}). Just as we want, these terms are the same such that we can use the cancelation of these leading terms in (\ref{3.5})-(\ref{3.8}) containing the wave operator disappears. In other words, we actually get the Huygens' wave with the desired extra decay rate for the difference of the second column and the forth column in the Green's function, which is crucial for us to deal with the nonlinear coupling in Section 4.

\begin{lemma}\label{l 3.2} For any $|\alpha|\geq0$, there exists a constant $C>0$ such that
\begin{equation*}
\begin{array}{rl}
&|D_x^\alpha(G_{12}^l-G_{14}^l,\ G_{32}^l-G_{34}^l)|
\leq C\underline{(1\!+\!t)^{-\frac{5+|\alpha|}{2}}\big(1+\frac{(|x|-ct)^2}{1+t}\big)^{-N}},\\[2mm]
&|D_x^\alpha(G_{22}^l-G_{24}^l,G_{42}^l-G_{44}^l)|
\leq C(1\!+\!t)^{-\frac{3+|\alpha|}{2}}\big(1+\frac{|x|^2}{1+t}\big)^{-\frac{3+|\alpha|}{2}}
\!+C\underline{(1\!+\!t)^{-\frac{5+|\alpha|}{2}}\big(1+\frac{(|x|-ct)^2}{1+t}\big)^{-N}},
\end{array}
\end{equation*}
where the constant $N>0$ can be arbitrarily large.
\end{lemma}
We emphasize that the underlined two terms above will be used to deal with the nonlinear coupling for the non-conservative nonlinear terms.

\vspace{3mm}
\textbf{Pointwise description of high frequency part.}\vspace{3mm}

From (\ref{2.27})-(\ref{2.28}), the high frequency part of the first and the third rows of the Green's function have the following estimates:
\begin{equation}\label{3.9}
\begin{array}{rl}
&\hat{G}_{11}^h\sim \frac{\mathrm{e}^{r_{1} t}}{|\xi|^2}+\mathrm{e}^{r_{3} t}+\mathrm{e}^{r_{4} t}+\frac{\mathrm{e}^{r_{2} t}}{|\xi|^5},\ \ \
\hat{G}_{12}^h\sim \Big(\frac{\mathrm{e}^{r_{1} t}}{|\xi|^3}+\frac{\mathrm{e}^{r_{2} t}}{|\xi|^5}+\mathrm{e}^{r_{3} t}+\mathrm{e}^{r_{4} t}\Big)\frac{\xi }{|\xi|},\\[2mm]
&\hat{G}_{13}^h\sim \frac{\mathrm{e}^{r_{1} t}+\mathrm{e}^{r_{3} t}+\mathrm{e}^{r_{4} t}}{|\xi|^2}+\frac{\mathrm{e}^{r_{2} t}}{|\xi|^4},\ \
\hat{G}_{14}^h\sim \Big(\frac{\mathrm{e}^{r_{1} t}+\mathrm{e}^{r_{2} t}}{|\xi|^3}+\frac{\mathrm{e}^{r_{3} t}+\mathrm{e}^{r_{4} t}}{|\xi|^2}\Big)\frac{\xi }{|\xi|},
\end{array}
\end{equation}
and
\begin{equation}\label{3.10}
\begin{array}{rl}
&\hat{G}_{31}^h\sim \frac{\mathrm{e}^{r_{1} t}+\mathrm{e}^{r_{3} t}+\mathrm{e}^{r_{4} t}}{|\xi|^2}+\frac{\mathrm{e}^{r_{3} t}}{|\xi|^4},\ \ \
\hat{G}_{32}^h\sim \Big(\frac{\mathrm{e}^{r_{1} t}+\mathrm{e}^{r_{2} t}}{|\xi|^3}+\frac{\mathrm{e}^{r_{2} t}+\mathrm{e}^{r_{4} t}}{|\xi|^4}\Big)\frac{\xi }{|\xi|},\\[2mm]
&\hat{G}_{33}^h\sim \mathrm{e}^{r_{1} t}+\frac{\mathrm{e}^{r_{2} t}}{|\xi|^2},\ \
\hat{G}_{34}^h\sim \Big(\frac{\mathrm{e}^{r_{1} t}+\mathrm{e}^{r_{2} t}}{|\xi|}+\frac{\mathrm{e}^{r_{3} t}+\mathrm{e}^{r_{4} t}}{|\xi|^2}\Big)\frac{\xi }{|\xi|}.
\end{array}
\end{equation}
Then from (\ref{3.3})-(\ref{3.4}), the rest two rows satisfy
\begin{equation}\label{3.11}
\begin{array}{rl}
&\hat{G}_{21}^h\sim \Big(\frac{\mathrm{e}^{r_{1} t}}{|\xi|^3}+\frac{\mathrm{e}^{r_{2} t}}{|\xi|^5}+\mathrm{e}^{r_{3} t}+\mathrm{e}^{r_{4} t}\Big)\frac{\xi }{|\xi|},\\[2mm]
&\hat{G}_{22}^h\sim \Big(\frac{\mathrm{e}^{r_{1} t}}{|\xi|^2}+\frac{\mathrm{e}^{r_{2} t}}{|\xi|^4}+\mathrm{e}^{r_{3} t}+\mathrm{e}^{r_{4} t}+\mathrm{e}^{\kappa_{1} t}\Big)\frac{\xi\xi^T }{|\xi|^2}+e^{-\kappa_1t}I,\\[2mm]
&\hat{G}_{23}^h\sim \Big(\frac{\mathrm{e}^{r_{1} t}+\mathrm{e}^{r_{2} t}}{|\xi|^3}+\frac{\mathrm{e}^{r_{3} t}+\mathrm{e}^{r_{4} t}}{|\xi|^2}\Big)\frac{\xi}{|\xi|},\\[2mm]
&\hat{G}_{24}^h\sim \Big(\frac{\mathrm{e}^{r_{1} t}}{|\xi|^4}+\frac{\mathrm{e}^{r_{2} t}+\mathrm{e}^{r_{3} t}+\mathrm{e}^{r_{4} t}+\mathrm{e}^{\kappa_{1} t}+\mathrm{e}^{\kappa_{2} t}}{|\xi|^2}\Big)\frac{\xi\xi^T }{|\xi|^2}+\frac{\mathrm{e}^{\kappa_{1} t}+\mathrm{e}^{\kappa_{2} t}}{|\xi|^2}I,
\end{array}
\end{equation}
and
\begin{equation}\label{3.12}
\begin{array}{rl}
&\hat{G}_{41}^h\sim \Big(\frac{\mathrm{e}^{r_{1} t}+\mathrm{e}^{r_{2} t}}{|\xi|^3}+\frac{\mathrm{e}^{r_{3} t}+\mathrm{e}^{r_{4} t}}{|\xi|^2}\Big)\frac{\xi}{|\xi|},\\[2mm]
&\hat{G}_{42}^h\sim \Big(\frac{\mathrm{e}^{r_{1} t}}{|\xi|^4}+\frac{\mathrm{e}^{r_{2} t}}{|\xi|^2}+\frac{\mathrm{e}^{r_{3} t}+\mathrm{e}^{r_{4} t}+\mathrm{e}^{\kappa_{1} t}+\mathrm{e}^{\kappa_{2} t}}{|\xi|^2}\Big)\frac{\xi\xi^T }{|\xi|^3}+\frac{\mathrm{e}^{\kappa_{1} t}+\mathrm{e}^{\kappa_{2} t}}{|\xi|^2}I,\\[3mm]
&\hat{G}_{43}^h\sim \frac{\mathrm{e}^{r_{1} t}+\mathrm{e}^{r_{2} t}+\mathrm{e}^{r_{3} t}+\mathrm{e}^{r_{4} t}}{|\xi|^2}\xi,\\[3mm]
&\hat{G}_{44}^h\sim \Big(\frac{\mathrm{e}^{r_{1} t}+\mathrm{e}^{r_{3} t}+\mathrm{e}^{r_{4} t}}{|\xi|^2}\Big)\frac{\xi\xi^T }{|\xi|^2}+\Big(\mathrm{e}^{r_{2} t}+\mathrm{e}^{\kappa_{2} t}\Big)\frac{\xi\xi^T }{|\xi|^2}+e^{\kappa_2t}I.
\end{array}
\end{equation}
With (\ref{3.9})-(\ref{3.12}) in hand, we can have the pointwise description of the high frequency part from the spectrum analysis in (\ref{2.24(1)}) and (\ref{2.38}), we can conclude that the high frequency part contains two kinds of singularity. One is like the heat kernel with the singularity at $t=0$, and the other is like $|\xi|^{\beta}e^{-t}$ with an integer $\beta\leq0$. Consequently, by using Lemma \ref{A.3}, we can get the pointwise description for the high frequency part as follows.
\begin{lemma}\label{l 3.4}
There exists a constant $C>0$ such that the high frequency part satisfies
\begin{equation}
\left| D_{x}^{\alpha}(G_{ij}^h(x,t)-S_{ij}(x,t))\right|\leq Ce^{-t/C}(1+|x|^2)^{-N},
\end{equation}
for any integer $N>0$. Here the singular part $S_{ij}$ satisfies
\begin{equation}
\begin{array}{rl}
&\D for\ \ (i,j)=(1,2),(2,1),\ \ S_{ij}(x,t)=Ce^{-t/C}\Big[t^{-\frac{3+|\alpha|}{2}}e^{-\frac{|x|^2}{Ct}}+\underline{D_x(\delta(x)+f(x))}\Big],\\[2mm]
&\D for\ \ (i,j)=(2,2),\ \ \ \ \ \ \ S_{ij}(x,t)=Ce^{-t/C}\Big[t^{-\frac{3+|\alpha|}{2}}e^{-\frac{|x|^2}{Ct}}+\underline{\delta(x)+D_x(\delta(x)+f(x))}\Big],\ \ \\[2mm]
&\D for\ \ (i,j)\neq(1,2),(2,1),(2,2),\ \  S_{ij}(x,t)=Ce^{-t/C}\Big[t^{-\frac{3+|\alpha|}{2}}e^{-\frac{|x|^2}{Ct}}+\underline{\delta(x)+f(x)}\Big],\\[2mm]
&\ and\ f(x)\in L^1,\ {\rm supp}{f}(x)\subset\{|x|\leq\eta_0\ll1\}.
\end{array}
\end{equation}
\end{lemma}

The middle frequency part of the Green's function are bounded and analytic, since the only possible pole has been excluded here. In other words, it does not matter the pointwise estimate of the Green's function. Hence,  we have the following
pointwise descriptions for the Green's function.

\begin{theorem}\label{l 3.5} For any $|\alpha|\geq0$, there exists a constant $C>0$ such that
\begin{equation}\label{3.38}
\begin{array}{rl}
&\ \ \ \ \ |D_x^\alpha(G_{11}\!-\!S_{11},G_{13}\!-\!S_{13},G_{24}\!-\!S_{24}, G_{31}\!-\!S_{31},G_{33}\!-\!S_{33},G_{44}\!-\!S_{44})|\\[2mm]
&\ \ \ \ \ \ \ \ \leq C(1+t)^{-\frac{3+|\alpha|}{2}}\Big(1+\frac{|x|^2}{1+t}\Big)^{\!-N}+C(1+t)^{-\frac{4+|\alpha|}{2}}\Big(1+\frac{(|x|-ct)^2}{1+t}\Big)^{\!-N},\\[2mm]
&|D_x^\alpha (G_{22}\!-\!S_{22},G_{42}\!-\!S_{42})|\leq C(1\!+\!t)^{-\frac{3+|\alpha|}{2}}\!\Big(1\!+\!\frac{|x|^2}{1+t}\Big)^{\!\!-\frac{3+|\alpha|}{2}}\!\!\!\!
+\!C(1\!+\!t)^{-\frac{4+|\alpha|}{2}}\!\Big(1\!+\!\frac{(|x|-ct)^2}{1+t}\Big)^{\!\!-N}\!\!,\\[3mm]
&|D_x^\alpha(G_{12}\!-\!S_{12},G_{14}\!-\!S_{14},G_{21}\!-\!S_{21},G_{23}\!-\!S_{23}, G_{32}\!-\!S_{32},G_{34}\!-\!S_{34},G_{41}\!-\!S_{41},G_{43}\!-\!S_{43})|\\[2mm]
&\ \ \ \ \ \leq C(1+t)^{-\frac{4+|\alpha|}{2}}\Big(1+\frac{(|x|-ct)^2}{1+t}\Big)^{\!-N},
\end{array}
\end{equation}
and
\begin{equation}\label{3.39}
\begin{array}{lll}
&\!\!\!\!\!\!\!\!\!\!\!|D_x^\alpha(G_{12}\!-\!G_{14}\!-\!(S_{12}\!-\!S_{14}),G_{32}\!-\!G_{34}\!-\!(S_{32}\!-\!S_{34}))|
\leq C(1\!+\!t)^{-\frac{5+|\alpha|}{2}}\Big(1\!+\!\frac{(|x|-ct)^2}{1+t}\Big)^{\!-N},\\[2mm]
&\!\!\!\!\!\!\!\!\!\!\!\!\!\!|D_x^\alpha(G_{22}\!-\!G_{24}\!-\!(S_{22}\!-\!S_{24}),G_{42}\!-\!G_{44}\!-\!(S_{42}\!-\!S_{44}))|\\[2mm]
&\ \ \ \ \leq C(1\!+\!t)^{-\frac{3+|\alpha|}{2}}\!\Big(1+\frac{|x|^2}{1+t}\Big)^{\!-\frac{3+|\alpha|}{2}}
\!+\!C(1+t)^{-\frac{5+|\alpha|}{2}}\Big(1\!+\!\frac{(|x|-ct)^2}{1+t}\Big)^{\!-N},
\end{array}
\end{equation}
where the positive constant $N$ can be arbitrarily large, and $S_{ij}$ is defined in Lemma \ref{3.4}.
\end{theorem}

\section{Pointwise estimate for nonlinear problem}
In this section, we consider the nonlinear problem. First of all, by using Duhamel's principle, we can get the representation of the solution $(\rho,m,n,\omega)$ for the nonlinear problem (\ref{2.1})-(\ref{2.2}):
\begin{equation}\label{4.1}
 D_x^\alpha\left(\!
            \begin{array}{c}
              \rho \\
              m \\
              n\\
              \omega
            \end{array}
          \!\right)\!=D_x^\alpha G\ast_x\!\left(
            \begin{array}{ccc}
              \rho_{0} \\
              m_{0}\\
              n_0\\
              \omega_{0}
            \end{array}
          \right)+\int_0^t\!D_x^\alpha G(\cdot,t\!-\!s)\!\ast_x\!\left(\!
                                            \begin{array}{ccc}
                                              0 \\
                                              F_1 \\
                                              0 \\
                                              F_2
                                            \end{array}
                                          \!\right)\!(\cdot,s)ds,
\end{equation}
where the nonlinear terms $F_1,F_2$ are defined in (\ref{2.2}).\newpage

\textbf{Initial propagation.}\ Theorem \ref{l 3.5} and the initial condition (\ref{1.3}) together with the representation (\ref{3.1})-(\ref{3.2}) yield the linear estimate as follows:
\begin{equation}\label{4.2}
\begin{array}{rl}
&\left|D_x^\alpha G\ast_x\!\!\left(\!
            \begin{array}{ccc}
              \rho_{0} \\
              m_{0}\\
              n_0\\
              \omega_{0}
            \end{array}
          \!\right)\right|
          \leq\!  \left|D_x^\alpha(G-G_S)
                           \!\ast_x\!\left(\!\!
            \begin{array}{ccc}
              \rho_{0} \\
              m_{0}\\
              n_0\\
              \omega_{0}
            \end{array}
          \!\right)\right|
          +\left| G_S\ast_x\!D_x^\alpha \left(\!
            \begin{array}{ccc}
              \rho_{0} \\
              m_{0}\\
              \omega_{0} \\
              z_0
            \end{array}
          \!\right)\right|\\
          \leq &\displaystyle C\epsilon\Big((1+t)^{-\frac{3+|\alpha|}{2}}\Big(1+\frac{|x|^2}{1+t}\Big)^{-\frac{3}{2}}
          +(1+t)^{-\frac{4+|\alpha|}{2}}\Big(1+\frac{(|x|-ct)^2}{1+t}\Big)^{-\frac{3}{2}}\Big),\ |\alpha|\leq1.
\end{array}
\end{equation}
Here we have used the convolution estimate in Lemma \ref{A.4} for the initial propagation. Recall that the singular parts $S_{12},S_{21},S_{22}$ are like $D_x\delta(x)$ and the others are basically like $\delta(x)$. Thus, when dealing with the convolution between these terms and the initial data, one should put the derivatives on the initial data, which is the reason why giving the different initial assumptions for $(\rho,m)$ and $(n,\omega)$ in (\ref{1.3}).

\vspace{3mm}
\textbf{Nonlinear coupling.}\ According to the initial propagation, we give the following ansatz for the nonlinear problem for $|\alpha|\leq 1$:
\begin{equation}\label{4.3}
|D_x^\alpha(\rho,m,n,\omega)|\leq 2C\epsilon\Big((1+t)^{-\frac{3+|\alpha|}{2}}\big(1+\frac{|x|^2}{1+t}\big)^{-\frac{3}{2}}
          +(1+t)^{-\frac{4+|\alpha|}{2}}\big(1+\frac{(|x|-ct)^2}{1+t}\big)^{-\frac{3}{2}}\Big).
\end{equation}
Here we need the pointwise ansatz for the solution and its first derivatives due to the quasi-linearity of the system together with the singularity in the high frequency part of the Green's function. Additionally, the $H^5$-regularity assumption on the initial data is also because we will meet the $L^\infty$-norm of the third derivatives of the unknowns when dealing with the convolution between the singular part of the Green's function and the nonlinear terms.

For the nonlinear coupling, we take the momentum $m$ for example, since its corresponding components $G_{21}$ and $G_{22}$ in the Green's function contains the highest order singularity (high frequency part). For simplicity, the nonlinear coupling of $m$ is denoted as $\tilde{m}$ and have
\begin{equation}\label{4.4}
D_x^\alpha\tilde{m}=\int_0^t D_x^\alpha G_{22}(\cdot,t-s)\ast F_1(\cdot,s)+D_x^\alpha G_{24}(\cdot,t-s)\ast F_2(\cdot,s)ds,\ |\alpha|\leq1.
\end{equation}
Notice in (\ref{2.1})-(\ref{2.2})that the last term in $F_1$ and the last term in $F_2$ have the opposite sign. Then we can write the convolution  into the regular part and the singular part as follows:
\begin{equation}\label{4.5}
\begin{array}{rl}
&G_{22}\ast_x F_1+G_{24}\ast_x F_2
=\underbrace{(G_{22}-S_{22})\ast F_1+(G_{24}-S_{24})\ast F_2}_{regular\ part}+\underbrace{S_{22}\ast F_1+S_{24}\ast F_2}_{singular\ part},\\[2mm]
=&(G_{22}-S_{22})\ast \Big(-\operatorname{div}\big(\frac{m\otimes m}{\rho+\bar{\rho}}\big)\Big)\\[2mm]
&+(G_{24}-S_{24})\ast \Big(-\operatorname{div}\big(\frac{w\otimes m}{n+\bar{n}}\big)-\bar{\mu}\Delta\frac{nw}{n+\bar{n}}-(\bar{\mu}+\bar{\lambda})\nabla\operatorname{div}\frac{nw}{n+\bar{n}}
-\nabla\big(\rho(n+\bar{n})-\alpha_{1}n\big)\Big)\\[2mm]
&+\big(G_{22}-G_{24}-(S_{22}-S_{24})\big)\ast \Big(\big(\frac{\rho+\bar{\rho}}{n+\bar{n}}-\frac{\bar{\rho}}{\bar{n}}\big)\omega\Big)\\[2mm]
&+S_{22}\ast F_1+S_{24}\ast F_2.
\end{array}
\end{equation}
After inserting the first two convolutions of (\ref{4.5}) into (\ref{4.4}), the resulting integral can be estimated by putting one derivative of the nonlinear terms onto the Green's function and using $K_3$ in Lemma \ref{A.5}. The nonlinear term in the third convolution of (\ref{4.5}) has not the divergence form, however, according to the cancelation of $G_{22}$ and $G_{24}$ in the low frequency part established in Lemma \ref{l 3.2}. Then it also can be estimated by using $K_3$ in Lemma \ref{A.5}. The first derivative of the convolution for the regular part can be treated similarly.

Next, we consider the singular part in the above convolution. Note first that $G_{22}\sim e^{-t}D_x\delta(x)$ from Lemma \ref{l 3.4}, then for $|\alpha|=1$,
\begin{equation}\label{4.6}
\begin{array}{rl}
&\D\Big|\int_0^t  S_{22}(\cdot,t-s)\ast D_xF_1(\cdot,s)ds\Big|
\sim \Big|\int_0^t \mathrm{e}^{-(t-s)}D_x\delta(\cdot,t-s)\ast_x D_x^2(|m|^2)(\cdot,s)ds\Big|\\
\leq &\D C\int_0^t\mathrm{e}^{-(t-s)}|D_x^3(|m|^2)|(x,s)ds\\
\leq &\D C\Big((1+t)^{-2}\Big(1+\frac{|x|^2}{1+t}\Big)^{-\frac{3}{2}}
          +(1+t)^{-\frac{5}{2}}\Big(1+\frac{(|x|-ct)^2}{1+t}\Big)^{-\frac{3}{2}}\Big),
\end{array}
\end{equation}
where we have used the decay rate of $\|D^3m\|_{L^\infty}$ and the ansatz (\ref{4.3}). In fact, the $L^\infty$-decay rate of $\|D^3m\|_{L^\infty}$ can be derived by using the $L^2$-decay rate of the solution and all of its derivatives in \cite{Wugc} together with the Sobolev inequality. The last convolution can be estimated similarly.

Up to now, we have closed the ansatz (\ref{4.3}) in $H^5$-framework, and get the pointwise space-time description of the solution and its first derivative. This proves Theorem \ref{l 1}.

\section {Appendix}

Some useful lemmas are given here. The first one describes the singular part of the high frequency:



\begin{lemma}\label{A.3}[Wang-Yang\cite{Wang}]
If ${\rm supp}\hat{f}(\xi)\subset
O_K=:{\{\xi, |\xi|\geq K>0\}}$, and $\hat{f}(\xi)$ satisfies
\begin{equation*}
|D_\xi^\beta\hat{f}(\xi)|\leq C|\xi|^{-|\beta|-1}\ \
({\rm or}\ |D_\xi^\beta\hat{f}(\xi)|\leq C|\xi|^{-|\beta|}),
\end{equation*}
then there exist distributions $f_1(x), f_2(x)$ and a constant $C_0$
such that $$
f(x)=f_1(x)+f_{2}(x)+C_{0}\delta(x)\ \
({\rm or}\ f(x)=f_1(x)+f_{2}(x)+C_{0}D_x\delta(x)),
$$
where $\delta(x)$ is the Dirac function. Furthermore, for any $|\alpha|\geq0$ and any positive integer $N$, we have
\begin{equation*}
|D_x^\alpha f_1(x)|\leq C(1+|x|^2)^{-N},\
\|f_{2}\|_{L^1}\leq C,\ {\rm supp}f_{2}(x)\subset\{x;|x|<\eta_0\ll1\}.
\end{equation*}

\end{lemma}

The last two lemmas are used for initial propagation and nonlinear coupling, respectively. We just state several typical cases here.

\begin{lemma}\label{A.4}[Wu-Wang\cite{Wu4}] There exists a constant $C>0$ such that:
\begin{equation*}
\begin{array}{ll}\label{6.6}
\displaystyle \int_{\mathbb R^3}\big(1\!+\!\frac{|x-y|^2}{1\!+\!t}\big)^{-n_1}\big(1+|y|^2\big)^{-n_2}dy\leq C\Big(1\!+\!\frac{|x|^2}{1+t}\Big)^{-n_3},\ \ {\rm for}\ n_1,n_2>\frac{3}{2}\ {\rm and}\ n_3=\min\{n_1,n_2\},\\
\displaystyle \int_{\mathbb R^3}\big(1+\frac{(|x-y|-ct)^2}{1\!+\!t}\big)^{-N}\big(1+|y|^2\big)^{-r_1}dy\leq C\Big(1+\frac{(|x|-ct)^2}{1\!+\!t}\Big)^{-\frac{3}{2}},\ \ {\rm for}\ N\geq r_1>\frac{21}{10}.
\end{array}
\end{equation*}
\end{lemma}

\begin{lemma}\label{A.5}[Liu-Noh\cite{LS}] There exists a constant $C>0$ such that
\begin{equation*}\label{6.1}
\begin{array}{ll}
\displaystyle K_1=\int_0^t\!\!\int_{\mathbb{R}^3}(1+t-s)^{-2}\Big(1+\frac{|x-y|^2}{1\!+\!t\!-\!s}\Big)^{-2}(1+s)^{-3}\Big(1+\frac{|y|^2}{1\!+\!s}\Big)^{-3}\!\!dyds
\leq C(1+t)^{-2}\big(1+\frac{|x|^2}{1+t}\big)^{-\frac{3}{2}},\\
\displaystyle K_2=\int_0^t\!\!\int_{\mathbb{R}^3}(1+t-s)^{-2}\Big(1+\frac{|x-y|^2}{1\!+\!t\!-\!s}\Big)^{-2}(1+s)^{-4}\Big(1+\frac{(|y|-cs)^2}{1\!+\!s}\Big)^{-3}dyds\\[1mm]
\ \ \ \ \leq \ C(1+t)^{-2}\Big(\big(1+\frac{|x|^2}{1+t}\big)^{-\frac{3}{2}}+\big(1+\frac{(|x|-ct)^2}{1+t}\big)^{-\frac{3}{2}}\Big),\\[1mm]
K_3=\displaystyle \int_{0}^{t}\!\!\int_{\mathbb{R}^3}(1+t-s)^{-\frac{5}{2}}\Big(1+\frac{(|x-y|-c(t-s))^2}{(1\!+\!t\!-\!s)}\Big)^{-N}(1+s)^{-4}\Big(1+\frac{(|y|-cs)^2}{1\!+\!s}\Big)^{-3}dyds\\[1mm]
\ \ \ \ \leq \ C(1+t)^{-2}\Big(\big(1+\frac{|x|^2}{1+t}\big)^{-\frac{3}{2}}+\big(1+\frac{(|x|-ct)^2}{1+t}\big)^{-\frac{3}{2}}\Big),
\end{array}
\end{equation*}
where the constant $N>0$ can be arbitrarily large.
\end{lemma}

\section*{Acknowledgments}

\bigbreak

{\bf Funding}: The research was supported by National Natural Science Foundation of China (No. 11971100) and the Fundamental Research Funds for the Central Universities (No. 2232019D3-43).\\
{\bf Conflict of Interest}: The authors declare that they have no conflict of interest.

\end{document}